\documentclass{amsart}
\usepackage{latexsym,amsfonts,amsmath,amssymb}
\newcommand{\Reals}{{\mathbb{R}}}
\newcommand{\Disk}{{\mathbb{D}}}
\newcommand{\Ints}{{\mathbb{Z}}}
\newcommand{\Cmplx}{{\mathbb{C}}}

\newcommand{\U}{{\mathbb{U}}}
\newcommand{\SO}{{\mathbb{SO}}}

\newcommand{\Exp}{{\mathbb{E}}}

\newcounter{smalllist}

\newtheorem{theorem}{Theorem}

\newtheorem{lemma}{Lemma}[section]
\newtheorem{prop}[lemma]{Proposition}
\newtheorem{coro}[lemma]{Corollary}
\theoremstyle{definition}
\newtheorem{definition}[lemma]{Definition}
\theoremstyle{remark}
\newtheorem*{remark}{Remark}

\DeclareMathOperator{\diag}{diag}
\DeclareMathOperator{\sign}{sign}

\let\llldots=\ldots
\def\ldots{\llldots{}}

\numberwithin{equation}{section}

\begin{document}

\title[Matrix models for circular ensembles]{Matrix models for circular ensembles}
\author[R.~Killip and I.~Nenciu]{Rowan Killip and Irina Nenciu}
\address{Rowan Killip\\
         UCLA Mathematics Department\\
         Box 951555\\
         Los Angeles, CA 90095}
\email{killip@math.ucla.edu}
\address{Irina Nenciu\\
         Mathematics 253-37\\
         Caltech\\
         Pasadena, CA 91125}
\email{nenciu@caltech.edu}
\date{\today}

\begin{abstract}
The Gibbs distribution for $n$ particles of the Coulomb gas on the unit
circle at
inverse temperature $\beta$ is given by\\[1mm]
\hbox{}\hfill\hbox{$\displaystyle
\mathbb{E}^{\beta}_{n}(f)=
\frac{1}{Z_{n,\beta}} \int\!\!\cdots\!\!\int\,
f(e^{i\theta_1},\ldots,e^{i\theta_n})|\Delta(e^{i\theta_1},\ldots,e^{i\theta_n}
)|^{\beta} \frac{d\theta_1}{2\pi}\cdots\frac{d\theta_n}{2\pi}
$}\hfill\hbox{}\\[1mm]
for any symmetric function $f$, where $\Delta$ denotes the Vandermonde
determinant and
$Z_{n,\beta}$ the normalization constant.

We will describe an ensemble of (sparse) random matrices whose eigenvalues follow this distribution.
Our approach combines elements from the theory of orthogonal polynomials on the
unit circle with ideas from recent work of Dumitriu and Edelman.
In particular, we resolve a question left open by them: find a tri-diagonal model for the Jacobi ensemble.
\end{abstract}

\maketitle

%%%%%%%%%%%%%%%%%%%%%%%%%%%%%%%%%%%%%%%%%%%%%%%%%%%%%%%%%%%%%%%%%%%%%%%%%%%%%%%%%%%%%
\section{Introduction}
%%%%%%%%%%%%%%%%%%%%%%%%%%%%%%%%%%%%%%%%%%%%%%%%%%%%%%%%%%%%%%%%%%%%%%%%%%%%%%%%%%%%%

In 1962, Dyson \cite{Dyson} introduced three ensembles of random unitary matrices with
a view to simplifying the study of energy level behavior in complex quantum systems.
Earlier work in this direction, pioneered by Wigner, focused on ensembles of Hermitian matrices.

The simplest of these three models is the unitary ensemble, which is just the group $\U(n)$
of $n\times n$ unitary matrices together with its Haar measure.  The induced probability
measure on the eigenvalues is given by the Weyl integration formula
(cf. \cite[\S VII.4]{Weyl}): for any symmetric function of the eigenvalues,
\begin{equation}\label{E:WeylIF}
  \Exp(f) = \frac{1}{n!} \int_0^{2\pi}\!\! \cdots \int_0^{2\pi}
f(e^{i\theta_1},\ldots,e^{i\theta_n})
    \bigl|\Delta(e^{i\theta_1},\ldots,e^{i\theta_n})\bigr|^2
    \frac{d\theta_1}{2\pi} \cdots \frac{d\theta_n}{2\pi}
\end{equation}
where $\Delta$ denotes the Vandermonde determinant,
\begin{equation}\label{VDefn}
  \Delta(z_1,\ldots,z_n) = \prod_{1\leq j < k \leq n} \!\! (z_k - z_j)
=
  \begin{vmatrix}
  1 & \ldots & 1 \\
  \vdots & & \vdots \\
  z_1^{n-1} & \cdots & z_n^{n-1}
  \end{vmatrix}.
\end{equation}

The orthogonal ensemble consists of symmetric $n\times n$ unitary matrices together with the unique measure that
is invariant under $U\mapsto W^T U W$ for all $W\in \U(n)$. Alternatively, if $U$ is chosen according
to the unitary ensemble, then $U^TU$ is distributed as a random element from the orthogonal ensemble.
The distribution of eigenvalues is given by \eqref{E:WeylIF} but with $|\Delta|^2$ replaced
by $|\Delta|$ and a new normalization constant.

The symplectic ensemble is a little more complicated.  Let $Z$ denote the $2n\times 2n$ block diagonal matrix
\begin{equation}
\begin{bmatrix}
  0  &  1 &      &     &    \\
 -1  &  0 &      &     &    \\
     &    &\ddots&     &    \\
     &    &      &  0  &  1 \\
     &    &      & -1  &  0
\end{bmatrix}
\end{equation}
and define the dual of a matrix by $U^R=Z^TU^TZ$.  The symplectic ensemble consists of self-dual unitary
$2n\times 2n$ matrices; the measure is that induced from $\U(2n)$ by the map $U\mapsto U^R U$.
(This is the unique measure invariant under $V\mapsto W^R V W$ for all $W\in \U(2n)$.)
The eigenvalues of such matrices are doubly degenerate and the pairs are distributed on the circle
as in \eqref{E:WeylIF} but now with $|\Delta|^2$ replaced by $|\Delta|^4$.  Again, the normalization
constant needs to be changed.

Dyson also observed that these eigenvalue distributions correspond to the Gibbs distribution for the classical
Coulomb gas on the circle at three different temperatures.  Let us elaborate.

Consider $n$ identically charged particles confined to move on the unit circle in the complex plane.
Each interacts with the others through the usual Coulomb potential, $-\log|z_i-z_j|$, which gives rise to
the Hamiltonian
$$
 H(z_1,\ldots,z_n) = \sum_{1\leq j<k\leq n}\!\! -\log|z_j-z_k|.
$$
(One may add a kinetic energy term; however, as we are interested only in the distribution of
the particle positions, it has no effect.)  This gives rise to the Gibbs measure (with parameters
$n$, the number of particles, and $\beta$, the inverse temperature)
\begin{align}
\Exp_n^\beta (f) &=
\frac{1}{(2\pi)^{n}Z_{n,\beta}} \int\!\! \cdots \!\! \int
f(e^{i\theta_1},\ldots,e^{i\theta_n})
    e^{-\beta H(e^{i\theta_1},\ldots,e^{i\theta_n})} \, d\theta_1\cdots
d\theta_n\\
       &= \frac{1}{(2\pi)^{n}Z_{n,\beta}} \int\!\! \cdots \!\! \int
f(e^{i\theta_1},\ldots,e^{i\theta_n}) \label{CGbeta}
    \bigl|\Delta(e^{i\theta_1},\ldots,e^{i\theta_n})\bigr|^\beta \,
d\theta_1 \cdots d\theta_n
\end{align}
for any symmetric function $f$. The partition function is given by
\begin{equation}\label{CGpart}
  Z_{n,\beta} = \frac{\Gamma(\tfrac12\beta n +
1)}{\bigl[\Gamma(\tfrac12\beta + 1)\bigr]^n}
\end{equation}
as conjectured by Dyson.  This was proved by Gunson \cite{Gunson} and Wilson \cite{Wilson},
though the Good proof \cite{Good} is even better.  We give another proof at the end of
Section~\ref{S:3}.

From the discussion above, we see that the orthogonal, unitary, and symplectic ensembles correspond to
the Coulomb gas at inverse temperatures $\beta=1$, $2$, and $4$.

Viewed from the opposite perspective, one may say that Dyson provided matrix models for the Coulomb gas
at three different temperatures.  Our first goal here is to present a family of matrix models
for all temperatures.  These matrices will be sparse---approximately $4n$ non-zero entries---which
suggests certain computational advantages.  To state the theorem, we need the following

\begin{definition}  We say that a complex random variable, $X$, with values in the unit disk,
$\Disk$, is $\Theta_\nu$-distributed (for $\nu>1$) if
\begin{equation}\label{E:ThetaDefn}
\Exp\{f(X)\} = \tfrac{\nu-1}{2\pi} \int\!\!\!\int_\Disk f(z)
(1-|z|^2)^{(\nu-3)/2} \,d^2z.
\end{equation}
For $\nu\geq2$ an integer, this has the following geometric interpretation: If $v$ is chosen from the
unit sphere $S^\nu$ in $\Reals^{\nu+1}$ at random according to the usual surface measure, then $v_1+iv_2$
is $\Theta_\nu$-distributed.  (See Corollary~\ref{C:A3}.)

As a continuation of this geometric picture, we shall say that $X$ is $\Theta_1$-distributed if it is
uniformly distributed on the unit circle in $\Cmplx$.
\end{definition}

Let us now describe the family of matrix models.

\begin{theorem}\label{T:1}
Given $\beta>0$, let $\alpha_k\sim\Theta_{\beta(n-k-1)+1}$ be independent random variables for $0\leq k\leq n-1$,
$\rho_k=\sqrt{1-|\alpha_k|^2}$, and define
$$
\Xi_k = \begin{bmatrix} \bar\alpha_k & \rho_k \\ \rho_k & -\alpha_k
\end{bmatrix}
$$
for $0\leq k\leq n-2$, while $\Xi_{-1}=[1]$ and $\Xi_{n-1}=[\bar\alpha_{n-1}]$ are $1 \times1$ matrices.
From these, form the $n\times n$  block-diagonal matrices
$$
L=\diag\bigl(\Xi_0   ,\Xi_2,\Xi_4,\ldots\bigr) \quad\text{and}\quad
M=\diag\bigl(\Xi_{-1},\Xi_1,\Xi_3,\ldots\bigr).
$$
Both $LM$ and $ML$ give {\rm(}sparse{\rm)} matrix models for the Coulomb gas at inverse temperature $\beta$.
That is, their eigenvalues are distributed according to \eqref{CGbeta}.
\end{theorem}

\begin{remark}
As each of the $\Xi_k$ is unitary, so are $L$ and $M$. (In the case of $\Xi_{n-1}$, we should reiterate
that $\alpha_{n-1}\sim\Theta_{1}$ is uniformly distributed on the unit circle.)  As a result, the eigenvalues
of $LM$ and $ML$ lie on the unit circle.  Note also that, since $M$ conjugates one to the other, $LM$ and
$ML$ have the same eigenvalues.
\end{remark}

In proving this theorem, we will be following the recent paper of Dumitriu
and Edelman \cite{DumE} rather closely, while incorporating the nuances
of the theory of polynomials orthogonal on the unit circle. The matrices $L$ and $M$ that appear in the theorem 
have their origin in the work of Cantero, Moral, and Vel\'asquez \cite{CMV}; this is discussed in Section~\ref{S:2}.

Dumitriu and Edelman constructed tri-diagonal matrix models for two of
the three standard examples of the Coulomb gas on the real line.
A model for the third will be constructed below.

The simplest way to obtain a normalizable Gibbs measure on the real line is to add an external harmonic
potential $V(x)=\tfrac{1}{2}x^2$.  This gives rise to the probability measure
\begin{equation}\label{CGRL}
\Exp(f) \propto
\int\!\! \cdots \!\! \int f(x_1,\ldots,x_n) \,
    \bigl|\Delta(x_1,\ldots,x_n)\bigr|^\beta  \prod_j e^{-V(x_j)}  \,
dx_1\cdots dx_n
\end{equation}
on $\Reals^n$.  This is known as the Hermite ensemble, because of its intimate connection to the
orthogonal polynomials of the same name, and when $\beta=1$, $2$, or $4$, arises as the eigenvalue
distribution in the classical Gaussian ensembles of random matrix theory.  Dumitriu and Edelman showed
that \eqref{CGRL} is the distribution of eigenvalues for a symmetric tri-diagonal matrix with independent
entries (modulo symmetry).  The diagonal entries have standard Gaussian distribution and the lower diagonal
entries are $2^{-1/2}$ times a $\chi$-distributed random variable with the number of degrees of freedom equal
to $\beta$ times the number of the row.

The second example treated by Dumitriu and Edelman is the Laguerre ensemble.  In statistical circles, this
is known as the Wishart ensemble, special cases of which arise in the empirical determination of the covariance
matrix of a multivariate Gaussian distribution. For this ensemble, one needs to modify the distribution given in
\eqref{CGRL} in two ways: each particle $x_j$ is confined to lie in $[0,\infty)$ and is subject
to the external potential $V(x)=-a\log(x)+x$, where $a>-1$ is a parameter.
In \cite{DumE}, it is shown that if $B$ is a certain $n\times n$ matrix with
independent $\chi$-distributed entries on the main and sub-diagonal (the number of degrees of freedom depends
on $a$, $\beta$, and the element in question) and zeros everywhere else, then the eigenvalues of
$L=BB^T$ follow this distribution.

The third canonical form of the Coulomb gas on $\Reals$ is the Jacobi ensemble.
The distribution is as in \eqref{CGRL}, but now the particles are
confined to lie within $[-2,2]$ and are subject to
the external potential $V(x)=-a\log(2-x)-b\log(2+x)$, where $a,b>-1$
are parameters.  This corresponds to the
probability measure on $[-2,2]^n$ that is proportional to
\begin{equation}\label{E:JE}
    \bigl|\Delta(x_1,\ldots,x_n)\bigr|^\beta  \prod_j \,
(2-x_j)^{a}(2+x_j)^{b}  \, dx_1\cdots dx_n.
\end{equation}
The partition function (or normalization coefficient) was determined by
Selberg~\cite{Selberg}.  This
will be discussed in Section~\ref{SandA}.

Dumitriu and Edelman did not give a matrix model for this ensemble, listing it as an open problem.
We present a tri-diagonal matrix model in Theorem~\ref{T:2} below.  The independent parameters follow
a beta distribution:

\begin{definition}
A real-valued random variable $X$ is said to be beta-distributed with parameters $s,t>0$,
which we denote by $X\sim B(s,t)$, if
\begin{equation}\label{E:beta}
\Exp\{f(X)\} = \frac{2^{1-s-t}\Gamma(s+t)}{\Gamma(s)\Gamma(t)}
\int_{-1}^1 f(x) (1-x)^{s-1}(1+x)^{t-1} \, dx.
\end{equation}

Note that $B(\tfrac{\nu}{2},\tfrac{\nu}{2})$ is the distribution of the
first component of a random vector from the $\nu$-sphere.  (See Corollary~\ref{C:A3}.)
\end{definition}

\begin{theorem}\label{T:2}
Given $\beta>0$, let $\alpha_k$, $0\leq k\leq 2n-2$, be distributed as follows
\begin{equation}\label{DistAR}
\alpha_k \sim \begin{cases} B(\tfrac{2n-k-2}{4}\beta + a +
1,\tfrac{2n-k-2}{4}\beta + b + 1)    &
\text{$k$ even,} \\
B(\tfrac{2n-k-3}{4}\beta + a+b+2,\tfrac{2n-k-1}{4}\beta)        & \text{$k$
odd.}
\end{cases}
\end{equation}
Let $\alpha_{2n-1}=\alpha_{-1}=-1$ and define
\begin{align}
b_{k+1}   &= (1-\alpha_{2k-1})\alpha_{2k} -
(1+\alpha_{2k-1})\alpha_{2k-2}          \label{E:BofG}\\
a_{k+1}   &= \big\{ (1-\alpha_{2k-1})(1-\alpha_{2k}^2)(1+\alpha_{2k+1})
\big\}^{1/2}    \label{E:AofG}
\end{align}
for $0\leq k \leq n-1$; then the eigenvalues of the tri-diagonal
matrix
$$
J = \begin{bmatrix}
 b_1 & a_1  &       &      \\
 a_1 & b_2  &\ddots &      \\
     &\ddots&\ddots & {\!a_{n-1}\!}     \\
     &      &{\!a_{n-1}\!}&  b_n
\end{bmatrix}
$$
are distributed according to the the Jacobi ensemble \eqref{E:JE}.
\end{theorem}

We know of two other papers which discuss the Jacobi ensemble in a manner inspired by
the work of Dumitriu--Edelman: \cite[\S 4.2]{FR} and \cite{Lip}.  These results are however of a rather
different character; in particular, we contend that Theorem~\ref{T:2} is the true Jacobi
ensemble analogue of the results of \cite{DumE}.

In Section~\ref{SandA}, we show how the ideas developed in the earlier parts of this paper
lead to new derivations of the classical integrals of Aomoto \cite{Aomoto} and Selberg \cite{Selberg}.
The main novelty of these proofs is their directness: They treat all values of $\beta$ and
$n$ on an equal footing.  In particular, we do not prove the result for $\beta$ an integer
and then make recourse to Carlson's Theorem.  These remarks are also applicable to the proof
of \eqref{CGpart} given at the end of Section~\ref{S:3}.

\smallskip

\noindent
\textit{Acknowledgements:} We would like to thank Barry Simon for encouragement and for
access to preliminary drafts of his forthcoming book \cite{Simon}.

%%%%%%%%%%%%%%%%%%%%%%%%%%%%%%%%%%%%%%%%%%%%%%%%%%%%%%%%%%%%%%%%%%%%%%%%%%%%%%%%%%%%%
\section{Overview of the Proofs and Background Material}\label{S:2}
%%%%%%%%%%%%%%%%%%%%%%%%%%%%%%%%%%%%%%%%%%%%%%%%%%%%%%%%%%%%%%%%%%%%%%%%%%%%%%%%%%%%%

We begin by examining the $\beta=2$ case of Theorem~\ref{T:1}, that is, Haar measure on the unitary group.

Rather than study the eigenvalues as the fundamental statistical object, we will consider the spectral
measure associated to $U$ and the vector $e_1=(1,0,\ldots,0)^T$. It will be denoted by $d\mu$.
As Haar measure is invariant under conjugation, any choice of unit vector $e_1$ leads to the same
probability distribution on $d\mu$.

The most natural coordinates for $d\mu$ are the eigenvalues $e^{i\theta_1},\ldots,e^{i\theta_n}$
and the mass that $d\mu$ gives to them: $\mu_1=\mu(\{e^{i\theta_1}\}),\ldots,\mu_{n-1}=\mu(\{e^{i\theta_{n-1}}\})$.
As $\int d\mu =1$, we omit $\mu_n=\mu(\{e^{i\theta_n}\})$.
Note that we have chosen not to order the eigenvalues, which means that the natural parameter space gives
an $n!$-fold cover of the set of measures.  We have already used this way of thinking a number of times,
beginning with \eqref{E:WeylIF}.

The above system of coordinates does not cover the possibility that $U$ has degenerate eigenvalues.  However,
as the Weyl integration formula shows, the set of such $U$ has zero Haar measure;
in fact, the density vanishes quadratically at these points.
The reason for this is worth repeating (cf.\ \cite[\S VII.4]{Weyl}): The submanifold where two eigenvalues
coincide has co-dimension three in $\U(n)$; one degree of freedom is lost from the reduction of the number
of eigenvalues and two more are lost in the reduction from two orthogonal one-dimensional eigenspaces to
a single two-dimensional eigenspace.  One should compare this to spherical polar coordinates in $\Reals^3$
where $r=0$ is a submanifold of co-dimension three and consequently, the density also vanishes to second order.

In Section~\ref{S:dmuUSO}, we will determine the probability distribution on $d\mu$ induced from Haar measure
on $\U(n)$, in the $\theta,\mu$ coordinates.  Conjugation invariance of Haar measure implies that
the eigenvalues and masses are statistically independent; it is then easy to see that the
former are distributed as in \eqref{E:WeylIF} and $(\mu_1,\ldots,\mu_n)$ is uniformly distributed on
the simplex $\sum \mu_j =1$.  See Proposition~\ref{P:Un}.  This implies that
$d\mu$ gives non-zero weight to each of the eigenvalues with probability one.  As a consequence, we can always
recover the eigenvalues from $d\mu$.

We will now introduce different coordinates, $(\alpha_0,\ldots,\alpha_{n-1})$, for $d\mu$
that arise in the study of orthogonal polynomials on the unit circle.

The monomials $1,z,\ldots,z^{n-1}$ form a basis for $L^2(d\mu)$ and so, applying the Gram--Schmidt procedure,
we can construct an orthogonal basis of monic polynomials: $\Phi_j$, $0\leq j<n$, with $\Phi_j$ monic of degree $j$.
We also define $\phi_j = \Phi_j / \|\Phi_j\|$, which gives an orthonormal basis.

There is a well-developed theory of such orthogonal polynomials, parts of which will be important in what follows.
For a proper discussion of this theory, see \cite{GeronBk}, \cite{Simon}, or \cite[Ch.~XI]{Szego}.

The first important fact about the orthogonal polynomials is that they
obey recurrence relations:
\begin{align}
\Phi_{k+1}(z)   &= z\Phi_k(z)   - \bar{\alpha}_k \Phi_k^*(z)
\label{PhiRec}\\
\Phi_{k+1}^*(z) &=  \Phi_k^*(z) -    \alpha_k z  \Phi_k(z)
\label{Phi*Rec}
\end{align}
where the $\alpha_k$ are the recurrence coefficients and $\Phi_k^*$
denotes the reversed polynomial:
\begin{equation}\label{rev}
\Phi_k(z) = \sum_{l=0}^k c_l z^l \quad \Rightarrow \quad \Phi_k^*(z) =
\sum_{l=0}^k \bar{c}_{k-l} z^l.
\end{equation}
Equivalently, $\Phi_k^*(z) = z^k \overline{\Phi_k(\bar{z}^{-1})}$.
These recurrence equations imply
\begin{equation}\label{PhiNorm}
\bigl\|\Phi_k\bigr\|_{L^2(d\mu)} = \prod_{l=0}^{k-1}
\rho_l\quad\text{where}\quad
    \rho_l=\sqrt{1-|\alpha_l|^2},
\end{equation}
from which the recurrence relations for the orthonormal polynomials are easily derived.

The recurrence coefficients $\alpha_k$ have been called by many names; we will follow \cite{Simon},
where they were recently dubbed `Verblunsky parameters'.  Each of $\alpha_0,\ldots,\alpha_{n-2}$ lies
inside the unit disk, $\Disk$, while $\alpha_{n-1}$ lies on its boundary, $S^1$.

There is an alternate way of relating measures to their Verblunsky parameters, namely the Schur algorithm:
If $d\mu$ is a probability measure, then we define its Schur function $f\!:\!\Disk\to\Disk$ by
\begin{equation}
f(z)=\frac{1}{z}\,\frac{F(z)-1}{F(z)+1} \quad\text{where} \quad
    F(z) = \int \frac{e^{i\theta}+z}{e^{i\theta}-z} d\mu(e^{i\theta}).
\end{equation}
The Schur algorithm parameterizes analytic maps $f\!:\!\Disk\to\bar\Disk$ by finitely or infinitely many
parameters $\alpha_k$.

There are finitely many parameters if and only if $f$ is a finite Blaschke product, or, equivalently,
iff $d\mu$ has finite support.  More precisely, the support of $d\mu$ consists of $n$ points
iff $f$ is a Blaschke product of degree $n$.  This is the case when there are $n$ Verblunsky parameters:
$\alpha_0,\ldots,\alpha_{n-2}\in \Disk$ and $\alpha_{n-1}\in S^1$. When there are finitely many
parameters, the last must always be unimodular.  In fact, the final parameter is essentially
equal to the product of the locations of the mass points of $d\mu$; see \eqref{detCMV}.

When $d\mu$ has infinite support, there are infinitely many Verblunsky
parameters, all of which lie in the unit disk.

Just as the Schur algorithm gives a bijection, so there is a bijection between measures $d\mu$ on $S^1$ supported
at $n$ points and sequences of parameters $\alpha_0,\ldots,\alpha_{n-2}\in \Disk$, $\alpha_{n-1} \in S^1$.  This
justifies their use as coordinates for the measure $d\mu$.

In Proposition~\ref{P:UnG}, we determine the probability distribution on $d\mu$ (induced by Haar measure on $\U(n)$)
in these new coordinates.  Interestingly, the $\alpha$s turn out to be statistically independent, with
$\alpha_k\sim\Theta_{2n-2k-1}$.

It is now but a few short steps to the $\beta=2$ case of Theorem~\ref{T:1}.

Consider the operator $f(z)\mapsto zf(z)$ in $L^2(d\mu)$.  The spectral measure associated to the
vector $f(z)\equiv 1$ is simply $d\mu$. To obtain a matrix model, we need only choose a basis
in which to represent this operator.  The most obvious choice is the basis of orthonormal polynomials, $\{\phi_k\}$.
This leads to a matrix whose entries can be expressed simply in terms of the $\alpha$s.
However, this matrix is not sparse: all entries above and including the sub-diagonal are non-zero (with probability one).  Such matrices are typically known as being in Hessenberg form.  In deference
to this, we will denote the matrix by $H$.  It plays an important role in the determination of the
distribution of the Verblunsky parameters, but does not appear in Theorem~\ref{T:1}.

The matrix $LM$ described in Theorem~\ref{T:1} is $f(z)\mapsto zf(z)$ in $L^2(d\mu)$ in the orthonormal basis
formed from $1,z,z^{-1},\ldots$ by applying the Gram--Schmidt procedure.
That this matrix can be expressed so simply in terms of the Verblunsky coefficients
is a discovery of Cantero, Moral, and Vel\'azquez, \cite{CMV}.
Related matters are discussed in Appendix~\ref{S:CMV}.  (The matrix $ML$ is the same operator in the basis
formed by applying Gram--Schmidt to $1,z^{-1},z,\ldots$.)

Thus far, we have only spoken about the unitary group, that is, about $\beta=2$.   In this case,
we have found a random ensemble of measures whose mass points are distributed as the particles in the
Coulomb gas at inverse temperature $\beta=2$.  The key discovery, however, was that the corresponding
Verblunsky parameters turn out to be independent.

For more general $\beta$, we wish to find an ensemble of measures so that the mass points are distributed
appropriately; we have complete freedom in choosing how the weights are distributed.  By the same token,
we want the induced probability distribution on the Verblunsky parameters to retain independence.  We
can then form the matrix set out in Theorem~\ref{T:1} and its eigenvalues are guaranteed to follow the
proper distribution.

The key to satisfying these desires is Lemma~\ref{VandUn}.  It expresses the value of
the Toeplitz determinant associated to $d\mu$ in terms of the $\theta,\mu$ coordinates and
in terms of the Verblunsky parameters.   Multiplying the probability distribution from the
$\beta=2$ case by the appropriate power of the Toeplitz determinant gives Proposition~\ref{ItU},
which is exactly the resolution of the goals set forth in the previous paragraph.

As an off-shoot of proving Theorem~\ref{T:1}, we are able to determine the Jacobian for the
map from the $(\theta,\mu)$-coordinates to the Verblunsky parameters $\alpha_k$.  That this is
possible is a delightful idea of Dumitriu and Edelman \cite{DumE}. (See Lemmas~\ref{L:Jacobian}
and~\ref{Somewhere}.)

Were we granted the Jacobian for this map, the paper could be much shorter---though we contend that the
scenic route followed below is not without merit.  It is a natural quantity to calculate and the
answer takes a rather simple form.  This behooves us to find a simple, direct derivation.  Thus far, we have
failed.  We would be much obliged to any reader who can resolve this matter.

The proof of Theorem~\ref{T:2} is very similar.  Again we begin by studying the problem for $\beta=2$.
The relevant group in this instance is not $\U(n)$, but rather $\SO(2n)$.  Such matrices have eigenvalues in
complex conjugate pairs and the corresponding eigenvectors are complex conjugates of one another.  Consequently,
the spectral measure associated to $e_1$ is symmetric with respect to complex conjugation.  The most
natural coordinates are $\theta_j\in (0,\pi)$ and $\mu_j\in[0,1]$ where
\begin{equation}\label{Ocoords}
    \int f d\mu = \sum_{j=1}^n \tfrac{1}{2}\mu_j[ f(e^{i\theta_j}) +
f(e^{-i\theta_j}) ]
\end{equation}
and $\sum \mu_j =1$.

Once again, we use the Verblunsky coefficients as a second set of
coordinates. These are now real as a consequence of
the complex conjugation symmetry of the measure.  Indeed, a measure has
this symmetry if and only if its Verblunsky
coefficients are real.  From this and the foregoing discussion of the
general case, we see that the last
Verblunsky coefficient, $\alpha_{2n-1}$, must be real and unimodular.
In fact it must be $-1$
because the product of the eigenvalues of a matrix from $\SO(2n)$ is
equal to one; see \eqref{detCMV}.
The remaining Verblunsky coefficients $\alpha_k$, $0\leq k \leq 2n-2$,
are free to range over $(-1,1)$.

By proceeding very much as before, we can construct certain ensembles of orthogonal matrices
for which the spectral measure is distributed in a desirable fashion.
When the Verblunsky coefficients are real, both $\det(1-U)$ and $\det(1+U)$
have simple expressions in terms of these coefficients.  As a result,
we are able to add two new parameters $a$ and $b$ to our family of distributions.
This line of reasoning leads us to Proposition~\ref{OIt}.

Given a measure, $d\mu$, on $S^1$ that is symmetric with respect to complex
conjugation, one may define a measure on $[-2,2]$ by
\begin{equation}\label{nuDefn}
  \int_{S^1} f(z+z^{-1}) \, d\mu(z)  =  \int_{-2}^2 f(x) \,d\nu(x).
\end{equation}
In particular, if $d\mu$ is of the form \eqref{Ocoords}, then we find
$$
\int f d\nu = \sum f(x_j) \mu_j \qquad\text{where $\quad x_j=2\cos(\theta_j)$.} 
$$
In this way, we find that Proposition~\ref{OIt} relates an
ensemble of probability measures on $[-2,2]$ to a certain ensemble of Verblunsky coefficients.
In fact, the locations of the masses of $d\nu$ are distributed as the points in the Jacobi
ensemble, \eqref{E:JE}, and are independent of the masses.

Theorem~\ref{T:2} follows immediately from the fact that the matrix $J$ represents
$f(x)\mapsto xf(x)$ in $L^2(d\nu)$ with respect to the basis of orthonormal polynomials.
The remainder of this section is devoted to explaining the origin of this fact.

Let $P_k(x)$ denote the monic polynomials orthogonal with respect to $d\nu$ and $p_k(x)$,
the corresponding orthonormal polynomials.  These obey a three-term recurrence relation \cite[\S3.2]{Szego}:
\begin{equation}\label{OPRLrec}
    x p_k(x) = a_{k+1} p_{k+1}(x) + b_{k+1} p_k(x) + a_{k} p_{k-1}(x),
\end{equation}
which explains the structure of the matrix $J$.

It is a famous observation of Szeg\H{o} (see \cite[\S11.5]{Szego}) that the polynomials orthogonal
with respect to $d\mu$ are intimately related to those orthogonal with respect to $d\nu$.
Specifically,
\begin{equation}\label{Sz}
  P_k(z+z^{-1}) = \frac{ z^{-k} \Phi_{2k}(z) + z^{k}\Phi_{2k}(z^{-1}) }{1-\alpha_{2k-1}}
\end{equation}
or equivalently,
$$
  p_k(z+z^{-1}) = \frac{z^{-k} \phi_{2k}(z) + z^{k} \phi_{2k}(z^{-1})}{\sqrt{2(1-\alpha_{2k-1})}}.
$$
Attendant to this relation between the orthogonal polynomials is a relation between their recurrence coefficients:
\begin{align*}
  b_{k+1}   &= (1-\alpha_{2k-1})\alpha_{2k} - (1+\alpha_{2k-1})\alpha_{2k-2}           \\
  a_{k+1}   &= \big\{ (1-\alpha_{2k-1})(1-\alpha_{2k}^2)(1+\alpha_{2k+1}) \big\}^{1/2}.
\end{align*}
These equations are known as the Geronimus relations. They are exactly equations
\eqref{E:BofG} and \eqref{E:AofG} from the statement of Theorem~\ref{T:2}.  This shows that the matrix $J$
defined in that theorem really does represent multiplication by $x$ in $L^2(d\nu)$.

The connection between the orthogonal polynomials for $d\mu$ and $d\nu$ suggests a relation between the matrices
$LM$ of Theorem~\ref{T:1} and $J$ of Theorem~\ref{T:2}.  This is investigated in Appendix~\ref{S:CMV}; in particular
we are able to give a short proof of the Geronimus relations.

%%%%%%%%%%%%%%%%%%%%%%%%%%%%%%%%%%%%%%%%%%%%%%%%%%%%%%%%%%%%%%%%%%%%%%%%%%%%%%%%%%%%%
\section{Distribution of $d\mu$ for $\U(n)$ and $\SO(2n)$}\label{S:dmuUSO}
%%%%%%%%%%%%%%%%%%%%%%%%%%%%%%%%%%%%%%%%%%%%%%%%%%%%%%%%%%%%%%%%%%%%%%%%%%%%%%%%%%%%%

Let $e_1$ denote the standard unit vector $(1,0,\ldots,0)^T$.  As described in the Introduction, for
each $U\in\U(n)$ we consider the spectral measure associated to the pair $(U,e_1)$, that is, the unique
measure on $S^1=\{z\in\Cmplx:|z|=1\}$ that obeys
$$
\langle e_1 | U^n e_1\rangle = \int z^n \,d\mu(z)
$$
for all $n\in\Ints$.  The first goal of this section is to determine the probability distribution of $d\mu$ when
$U\in\U(n)$ is chosen according to Haar measure. We will then prove the analogous result for $\SO(2n)$.
In each case, we give the distribution both in terms of the natural parameters
(the eigenvalues and associated masses) and in terms of the Verblunsky parameters.

\begin{prop}\label{P:Un}
If $U$ is chosen according to Haar measure on $\U(n)$, then the probability measure on $d\mu$ is given by
\begin{equation}\label{Uinzmu}
\frac{(n-1)!}{n! (2\pi)^n}
\bigl|\Delta(e^{i\theta_1},\ldots,e^{i\theta_n})\bigr|^2 d\theta_1 \cdots d\theta_n \, d\mu_1 \ldots d\mu_{n-1}
\end{equation}
where $0\leq \theta_j\leq 2\pi$ and $0\leq\mu_j\leq 1$ with $\sum \mu_j \leq 1$.
\end{prop}

\begin{proof}
Conjugation invariance of Haar measure on $\U(n)$ implies that the distribution of the
eigenvectors is independent of that of the eigenvalues.  The distribution of the eigenvalues is
given by the Weyl integration formula \eqref{E:WeylIF}, while the masses are simply the square moduli
of the top entries of the normalized eigenvectors.

Of course, the normalized eigenvectors are only determined up to a phase factor.  If we choose this phase
factor at random from the unit circle, then by conjugation invariance, the top entries of
the eigenvectors are distributed as a random unit vector from $\Cmplx^n$.  Therefore, by Corollary~\ref{C:A3},
$(\mu_1,\ldots,\mu_n)$ is uniformly distributed on the $(n-1)$-simplex $\sum \mu_j =1$.
\end{proof}

The key to writing the probability measure on $d\mu$ in the Verblunsky coordinates is the Householder
algorithm \cite[\S6.4]{Househ}, which converts any matrix to one in Hessenberg form (i.e., with zeros
below the sub-diagonal) via unitary conjugation.  The algorithm proceeds iteratively, killing the
undesirable elements in each column in turn, working left to right.

Consider a matrix $A$ whose first $k-1$ columns are zero below the sub-diagonal.
Let us write $A_k=[a_{1,k},\ldots,a_{n,k}]^T$ for the $k$th column of $A$ and then define
\begin{equation}\label{E:vDef}
v = [0,\ldots,0,\alpha,a_{k+2,k},\ldots,a_{n,k}]^T
\end{equation}
\begin{equation}\label{E:alphaDef}
    \alpha = a_{k+1,k} - \frac{a_{k+1,k}}{|a_{k+1,k}|}
\sqrt{|a_{k+1,k}|^2+\cdots+|a_{n,k}|^2}.
\end{equation}
The reflection through the plane perpendicular to $v$ is given by $R=I-2\frac{vv^\dagger}{\|v\|^2}$.
It maps $A_k$ to $A_k-v$, which has zeros in the desired places. Moreover, the first $k-1$ columns
of $A$ are unchanged by left multiplication by $R$ because they are orthogonal to $v$---indeed, their only
non-zero entries coincide with zeros in $v$.  Similarly, right multiplication of any matrix
by $R$ leaves its first $k$ columns unchanged.

As $R$ is a reflection, $R^{-1}=R=R^\dagger$ and so we find that conjugating $A$ by $R$ gives a new matrix
whose first $k$ columns agree with Hessenberg form (i.e., with zeros below the sub-diagonal).    In this
way, we have described one step of the usual Householder algorithm. However, we wish to add one further
conjugation so as to make the entries on the sub-diagonal non-negative. To do this, we form $DRARD^\dagger$
where $D$ differs from the identity matrix by having $(k+1,k+1)$ entry $e^{-i\phi}$ with $\phi$
chosen appropriately.

Starting with a unitary matrix $U$, we can apply the above algorithm to obtain a unitarily equivalent
matrix, $H$, in Hessenberg form with non-negative sub-diagonal. Moreover, the spectral measure for
$(H,e_1)$ is the same as that for $(U,e_1)$---namely $d\mu$---because the vector $e_1$ is fixed by all
the unitary matrices by which $U$ is conjugated.

In the Introduction, we used $H$ to denote the matrix representation of $f(z)\mapsto zf(z)$ on $L^2(d\mu)$
in the basis of orthonormal polynomials.  It is easily seen to be in Hessenberg form and to have a non-negative
(indeed positive) sub-diagonal.  As one might hope from the notation, these two matrices are the same (see also
\cite[Corollary 3.3]{CMV2}):

\begin{lemma}
Suppose $H$ is an $n\times n$ unitary matrix in Hessenberg form with non-negative sub-diagonal and let
$d\mu$ denote the spectral measure associated to the vector $e_1$.  If the support of $d\mu$ consists
of $n$ points, then $H$ represents $f(z)\mapsto zf(z)$ in the basis of orthonormal polynomials.  Consequently,
\begin{equation}\label{GGTentry}
H_{i+1,j+1} = \langle \phi_i | z \phi_j \rangle =
    \begin{cases}
    -\alpha_{i-1}\bar\alpha_{j} \prod_{l=i}^{j-1} \rho_{l} & i < j+1 \\
    \rho_{j-1} & i = j+1 \\
    0 & i > j+1
    \end{cases}
\end{equation}
where $\rho_j=\sqrt{1-|\alpha_j|^2}$ and $\alpha_{-1}=-1$.
\end{lemma}

\begin{proof}
As $d\mu$ is the spectral measure for $(H,e_1)$ and as $L^2(d\mu)$ has the same dimension as the space
on which the operator $H$ acts, there must be an orthonormal basis $f_1,\ldots,f_n$ for $L^2(d\mu)$
with $f_0\equiv 1$ such that $H$ represents $f(z)\mapsto zf(z)$ in this basis.  This is just the spectral
theorem combined with the fact that $e_1$ must be cyclic (for otherwise, $L^2(d\mu)$ wouldn't have full
dimension).

From the cyclicity argument we also learn that $H$ must have a strictly positive sub-diagonal.

To finish the proof of the first claim, we need only show that $f_j(z)=\phi_{j-1}(z)$; that is,
that the orthonormal basis in question is precisely that of the orthonormal polynomials.  Because $H$
is in Hessenberg form with positive sub-diagonal, the standard basis vectors arise from applying the
Gram--Schmidt procedure to $e_1,He_1,\ldots,H^{n-1}e_1$.  Consequently, the vectors $f_j$ must be
the result of applying the same procedure to $1,z,\ldots,z^{n-1}$; that is, $f_j$ must be $\phi_{j-1}$.

The first part of \eqref{GGTentry} merely reexpresses what we have just proved.  The second follows from
the recursion relations and \eqref{PhiNorm}; however, the proof is not particularly enlightening and
we refer the reader to \cite{GeronMS} for details.
\end{proof}

We will now apply the Householder algorithm outlined earlier to a matrix chosen at random from $\U(n)$.
By the lemma above, this will allow us to determine the induced distribution on the Verblunsky parameters
associated to the spectral measure for $(U,e_1)$.

\begin{prop}\label{P:UnG}
Let $d\mu$ be the spectral measure for $(U,e_1)$ with $U$ chosen at
random from $\U(n)$ according to Haar measure.
In terms of the Verblunsky parameters, $\alpha_0,\ldots,\alpha_{n-2}$,
and $\alpha_{n-1}=e^{i\phi}$, this probability
distribution is given by
\begin{equation}\label{Uingamma}
\tfrac{1}{2}\tfrac{(n-1)!}{\pi^n}\prod_{k=0}^{n-2}
(1-|\alpha_k|^2)^{n-2-k}
    \, d^2 \alpha_0\cdots d^2\alpha_{n-2} \, d\phi.
\end{equation}
That is, the Verblunsky parameters are independent and
$\alpha_j\sim\Theta_{2n-2j-1}$.
\end{prop}

\begin{proof}
The key to applying the Householder algorithm to a random element
$U\in\U(n)$ is the following realization of Haar measure:
Choose the first column at random from the unit sphere; then choose the
second column from the unit sphere of
vectors orthogonal to the first; then the third column and so forth.
In this way, one could say that the columns of $U$
form a random orthonormal basis for $\Cmplx^n$.  (That this is indeed
Haar measure is a simple consequence
of invariance under left multiplication by unitary matrices.)

The first column of $U$ is a random vector from the unit sphere.  After
left multiplication by the appropriate reflection $R$,
the new first column takes the form $[\bar\alpha_0,b,0,\ldots,0]^T$
where $\bar\alpha_0$ is the the original $(1,1)$ entry
of $U$ and so $\Theta_{2n-1}$-distributed, while $b$ has modulus
$\rho_0$ and arbitrary argument.  Subsequent left
multiplication by $D$ converts the first column to
$[\bar\alpha_0,\rho_0,0,\ldots,0]^T$, as it will remain.
The other columns are still orthogonal to the first; indeed,
they form a random orthonormal basis for the orthogonal complement of
the first column.

Right multiplication by $RD^\dagger$ leaves the first column untouched
while orthogonally intermixing the other columns.
Of course, this means that they remain a random orthonormal basis for
the orthogonal complement of the first column.  (Remember, Haar measure
is also invariant under right multiplication by a unitary.)

For the subsequent columns, the procedure is similar.  Let us skip
ahead to dealing with the $k$th column.

From the unitarity of the matrix $H$ from the previous lemma,
$$
X=
\begin{bmatrix}
   \rho_0\rho_1\rho_2\cdots\rho_{k-2} \\
-\alpha_0\rho_1\rho_2\cdots\rho_{k-2} \\
-\alpha_1\rho_2\cdots\rho_{k-2} \\
\vdots \\
-\alpha_{k-3}\rho_{k-2} \\
-\alpha_{k-2} \\
 0 \\
 \vdots \\
 0
\end{bmatrix}
$$
is a unit vector orthogonal to the first $k-1$ columns.  As the $k$th
column is a random vector orthogonal
to the first $k-1$ columns, its inner product with $X$ is distributed
as the top entry of a random vector
from the $(2n-2k+1)$-sphere and is independent of
$\alpha_1,\ldots,\alpha_{k-2}$.  Let us call this inner
product $\bar\alpha_{k-1}$, noting that this implies $\alpha_{k-1}$ is
$\Theta_{2n-2k+1}$-distributed
as stated in the proposition.

We now multiply the matrix at hand from the left by the appropriate
reflection and rotation to bring the $k$th
column into the desired
form.  Neither of these operations alters the top $k$ rows and so the
inner product of the $k$th column
with $X$ is unchanged.  But now the $k$th column is uniquely
determined; it must be $\bar\alpha_{k-1}X+\rho_{k-1}e_{k+1}$,
just as in \eqref{GGTentry}.

Lastly, we should multiply on the right by $RD^\dagger$, but this
leaves the first $k$ columns unchanged while
orthogonally intermixing the other columns.  In this way, we obtain a
matrix whose first $k$ columns conform to the
structure of $H$, while the remaining columns form a random basis for
the orthogonal complement of the span of
those $k$ columns.

In this way, we can proceed inductively until we reach the last column.
It is obliged to be a random orthonormal basis
for the one-dimensional space orthogonal to the preceding $n-1$ columns
and hence a random unimodular multiple,
say $\bar\alpha_{n-1}$, of $X$.  This is why the last Verblunsky
parameter is $\Theta_1$-distributed.

We have now conjugated $U$ to a matrix in the form of \eqref{GGTentry}
with parameters distributed as stated
in the proposition.  The vector $e_1$ is unchanged under the action of
each of the conjugating matrices and
consequently, these are precisely the Verblunsky parameters of $d\mu$.
\end{proof}

We now turn to the study of Haar measure on $\SO(2n)$.  The proofs
follow those given above pretty closely.

\begin{prop}\label{SO2n}
If $U$ is chosen at random from $\SO(2n)$ according to Haar measure,
then the spectral measure
$d\mu$ associated to $(U,e_1)$ is distributed as
\begin{equation}\label{SOinzmu}
\frac{(n-1)!}{2^{n-1}\ n!}
\bigl|\Delta\bigl(2\cos\theta_1,\ldots,2\cos\theta_n\bigr)\bigr|^2
    \frac{d\theta_1}{\pi}\cdots\frac{d\theta_n}{\pi} d\mu_1\cdots
d\mu_{n-1}
\end{equation}
where $\theta_j$ and $\mu_j$ are the coordinates given in
\eqref{Ocoords}.
\end{prop}

\begin{proof}
By the Weyl integration formula for $\SO(2n)$, the marginal distribution of the
eigenvalues is as above. (See \cite[\S VII.9]{Weyl}.)

If the eigenvalues are prescribed, say $e^{\pm i\theta_1},\ldots,e^{\pm i\theta_n}$, then
the conditional distribution of $U$ is given by taking a fixed matrix with this spectrum and
conjugating it by a random element from $\SO(2n)$.  The natural choice for this fixed matrix is
block diagonal:
$$
U_0 = \diag\left(
    \begin{bmatrix} \cos(\theta_1) & \sin(\theta_1) \\ -\sin(\theta_1) &
\cos(\theta_1)\end{bmatrix},\ldots,
    \begin{bmatrix} \cos(\theta_n) & \sin(\theta_n) \\
-\sin(\theta_n)&\cos(\theta_n)\end{bmatrix}\right).
$$
From this we see that the $d\mu$ is the spectral measure for $U_0$ and a random vector from the $(2n-1)$-sphere.
The proposition then follows by diagonalizing $U_0$ and applying Corollary~\ref{C:A3}.
\end{proof}

\begin{prop}\label{SO2nG}
Let $U$ be chosen from $\SO(2n)$ according to Haar measure and let
$d\mu$ denote the spectral measure
for $(U,e_1)$.  In terms of the Verblunsky parameters, the probability
distribution on $d\mu$ is
\begin{equation}\label{SOing}
\tfrac{(n-1)!}{\pi^n} \prod_{k=0}^{2n-2}
(1-\alpha_k^2)^{\frac{2n-k-3}{2}}\,d\alpha_0\cdots d\alpha_{2n-2}
\end{equation}
with $\alpha_{2n-1}=-1$.
That is, the Verblunsky parameters are independent and $\alpha_k\sim
B(\frac{2n-k-1}{2},\frac{2n-k-1}{2})$.
\end{prop}

\begin{proof}
While we may use the Householder algorithm as set out above, the fact
that we are now dealing with
real-valued matrices allows the following simplification: the vector
defining the reflection is
as in \eqref{E:vDef}, but now with
$$
\alpha = a_{k+1,k} - \sqrt{a_{k+1,k}^2 + \cdots + a_{2n,k}^2}
$$
instead of \eqref{E:alphaDef}.  This permits us to forgo the
conjugation by $D$.

Haar measure on $\SO(2n)$ can be realized by choosing the first column
as a random vector from the
unit sphere in $\Reals^{2n}$, and then the second as a random vector
orthogonal to the first, and so
forth.  However, the fact that the matrix has determinant one means that
the first $2n-1$ columns
completely determine the last.  One may say that the columns of $U$
form a random positively oriented
basis for $\Reals^{2n}$.

Proceeding as in the proof of Proposition~\ref{P:UnG}, we see that
$\alpha_{k-1}$ is defined as the
inner product of a specific vector $X$ with a random unit vector from
the $(2n-k-2)$-sphere of vectors
orthogonal to the first $k-1$ columns.  It follows from
Corollary~\ref{C:A2} that $\alpha_{k-1}\sim
B(\frac{2n-k-2}{2},\frac{2n-k-2}{2})$ as stated above.

The last column of $H$, and hence $\alpha_{2n-1}$, is uniquely determined by the fact that $\det(H)=1$.
It is just a matter of using \eqref{GGTentry} to determine which value of $\alpha_{2n-1}$ makes
this determinant one; moreover, by continuity of the determinant, it suffices to consider the case where
all other $\alpha$s are zero.  This gives
$$
1=\det(H)=-\alpha_{-1}\alpha_{2n-1}\sign(\sigma)=-\alpha_{2n-1}
$$
where $\sigma$ is the cyclic permutation $j\mapsto j+1 \mod 2n$, which is odd.

Lastly, we should justify the normalization coefficient given in \eqref{SOing}; what appears there
is very much simpler than one would expect from \eqref{E:beta}.  This
simplification is based on the duplication formula for the $\Gamma$ function:
$
\sqrt{\pi}\Gamma(2t) = 2^{2t-1}\Gamma(t)\Gamma(t+\frac12).
$
Specifically, beginning with \eqref{E:beta},
$$
\int_{-1}^1 (1-\alpha^2)^{t-1} \, d\alpha =
\frac{2^{1-2t}\Gamma(2t)}{\Gamma(t)^2}
    = \frac{\Gamma(t+\frac12)}{\sqrt{\pi}\,\Gamma(t)},
$$
which causes the product of normalization coefficients to telescope:
$$
      \prod_{k=0}^{2n-2} \frac{\Gamma(\frac{2n-k}{2}) }{ \sqrt{\pi}\,\Gamma(\frac{2n-k-1}{2}) }
    = \pi^{\frac{1}{2}-n} \frac{\Gamma(n)}{\Gamma(\tfrac{1}{2})}
    = \frac{(n-1)!}{\pi^n},
$$
as given in \eqref{SOing}.
\end{proof}

%%%%%%%%%%%%%%%%%%%%%%%%%%%%%%%%%%%%%%%%%%%%%%%%%%%%%%%%%%%%%%%%%%%%%%%%%%%%%%%%%%%%%%%%%%%%%%%%%%%%%%%%%
\section{The Proof of Theorem 1}\label{S:3}
%%%%%%%%%%%%%%%%%%%%%%%%%%%%%%%%%%%%%%%%%%%%%%%%%%%%%%%%%%%%%%%%%%%%%%%%%%%%%%%%%%%%%%%%%%%%%%%%%%%%%%%%%

Let $d\mu$ be the measure on $S^1$ given by
\begin{equation}\label{S3mu}
\int f d\mu = \sum_{j=1}^n \mu_j f(e^{i\theta_j})
\end{equation}
with $\theta_j\in[0,2\pi)$ distinct and $\sum \mu_j=1$.
As discussed in the Introduction, this measure is uniquely determined
by its Verblunsky parameters $\alpha_0,\ldots,\alpha_{n-2}\in\Disk$ and
$\alpha_{n-1}=e^{i\phi}\in S^1$.

It is difficult to find functions of $d\mu$ that admit simple expressions in terms of both
$\theta_j,\mu_j$ and the Verblunsky parameters.  One such quantity is the determinant of the
associated Toeplitz matrix; this is the subject of the next lemma.

\begin{lemma}\label{VandUn}
If $d\mu$ is a probability measure of the form given in \eqref{S3mu} and $\alpha_0,\ldots,\alpha_{n-1}$
its Verblunsky coefficients, then
\begin{equation}\label{VUnit}
|\Delta(z_1,\ldots,z_n)|^2 \, \prod_{j=1}^n \mu_j = \prod_{k=0}^{n-2}
(1-|\alpha_k|^2)^{n-k-1}.
\end{equation}
\end{lemma}

\begin{proof}
Let $c_k=\sum_{j=1}^{n} \mu_j z_j^k$ denote the moments of $d\mu$.
We will prove that both sides of \eqref{VUnit} are equal to the
determinant of the
$n\times n$ Toeplitz matrix associated to $d\mu$:
$$
T=\begin{bmatrix}
  c_0 & c_{-1} & \cdots & c_{1-n} \\
  c_1 & c_0 & \cdots & c_{2-n} \\
  \vdots & \vdots & \ddots & \vdots \\
  c_{n-1} & c_{n-2} & \cdots & c_0 \\
\end{bmatrix}.
$$

If we define $A$ and $M$ by
$$
A=\begin{bmatrix}
  1 & 1 & \cdots & 1 \\
  z_1 & z_2 & \cdots & z_n \\
  \vdots & \vdots & \ddots & \vdots \\
  z_1^{n-1} & z_2^{n-1} & \cdots & z_n^{n-1}
\end{bmatrix}
\quad
M=\begin{bmatrix}
  \mu_1 &   0   & \cdots & 0 \\
    0   & \mu_2 & \cdots & 0 \\
 \vdots &\vdots & \ddots & \vdots \\
    0   &   0   & \cdots & \mu_n
\end{bmatrix},
$$
then $T = A M A^\dagger$. Consequently,
$$
  \det(T) = \bigl|\det(A)\bigr|^2 \, \det(M) =
|\Delta(z_1,\ldots{},z_n)|^2 \, \prod_{j=1}^n \mu_j.
$$

We will now show that the right-hand side of \eqref{VUnit} is equal to
$\det(T)$.
To this end, let
$$
B=\begin{bmatrix}
  \Phi_0(z_1) & \Phi_0(z_2) & \cdots & \Phi_0(z_n) \\
  \Phi_1(z_1) & \Phi_1(z_2) & \cdots & \Phi_1(z_n) \\
  \vdots & \vdots & \ddots & \vdots \\
  \Phi_{n-1}(z_1) & \Phi_{n-1}(z_2) & \cdots & \Phi_{n-1}(z_n)
\end{bmatrix},
$$
which has the same determinant as $A$ because each can be reduced to the other by elementary row
operations.  From the orthogonality property of the $\{\Phi_j\}$, it follows that
$B M B^\dagger$ is the diagonal matrix whose entries are the squares of the $L^2(d\mu)$-norms of
$\Phi_0,\Phi_1,\ldots{,}\Phi_{n-1}$. Therefore by \eqref{PhiNorm},
$$
\det(T) = \det(B M B^\dagger) = \prod_{k=0}^{n-1} \bigl\|\Phi_k\bigr\|^2_{L^2(d\mu)}
= \prod_{k=0}^{n-2} (1-|\alpha_k|^2)^{n-k-1},
$$
just as was required.
\end{proof}

Both expressions for the Toeplitz determinant are well known; indeed, the argument presented above and its 
Hankel-matrix analog play a central role in random matrix theory.

We are now ready to state and prove the main result of this section.
Please note that neither measure given below is normalized; however, they do have the same
normalization coefficient.  It is calculated in Lemma~\ref{Somewhere} where it is used to
give an independent proof of \eqref{CGpart}.

\begin{prop}\label{ItU}
The following formulae express the same measure on the manifold of
probability distributions on $S^1$ supported at $n$ points:
\begin{equation}\label{qqqqq}
\frac{2^{1-n}}{n!} |\Delta(e^{i\theta_1},\ldots,e^{i\theta_n})|^{\beta}
\prod_{j=1}^n
\mu_j^{\frac{\beta}{2}-1} \, d\theta_1 \cdots d\theta_n\,d\mu_1 \cdots
d\mu_{n-1}
\end{equation}
in the $(\theta,\mu)$-coordinates of \eqref{S3mu}, and
\begin{equation}\label{ggggg}
\prod_{k=0}^{n-2} (1-|\alpha_k|^2)^{\frac{\beta}{2}(n-k-1)-1} \, d^2
\alpha_0 \cdots d^2\alpha_{n-2}\,d\phi
\end{equation}
in terms of the Verblunsky parameters.
\end{prop}

\begin{proof}
When $\beta=2$, this follows immediately from Propositions~\ref{P:Un} and~\ref{P:UnG}.
To obtain the general-$\beta$ version of \eqref{qqqqq} from the $\beta=2$ version, one
has to multiply by
\begin{equation}\label{VUnitLHS}
    |\Delta(e^{i\theta_1},\ldots,e^{i\theta_n})|^{\beta-2} \prod_{j=1}^n
\mu_j^{\frac{\beta}{2}-1},
\end{equation}
while the same transformation of $\beta$ in \eqref{ggggg} is effected by multiplying by
\begin{equation}\label{VUnitRHS}
    \prod_{k=0}^{n-2} (1-|\alpha_k|^2)^{(\frac{\beta}{2}-1)(n-k-1)}.
\end{equation}
But, \eqref{VUnitLHS} and \eqref{VUnitRHS} are equal; they are either side of \eqref{VUnit}
raised to the power $\frac{\beta}{2}-1$.
\end{proof}

\begin{proof}[Proof of Theorem 1]
Theorem~\ref{T:1} is an immediate corollary of Proposition~\ref{ItU} and results from \cite{CMV}:
The Verblunsky parameters of the spectral measure for $(LM,e_1)$ are precisely the $\alpha$s
that appear in the definition of $L$ and $M$.  Consequently, if the Verblunsky parameters
are distributed according to \eqref{ggggg}, then the eigenvalues are distributed as in \eqref{CGbeta}.
\end{proof}

It is fair to suggest that studying the unitary group is a rather roundabout way of proving the above proposition.
We simply do not know a better way. The natural suggestion is to first calculate the Jacobian for the map
from the $(\theta,\mu)$-coordinates to the Verblunsky parameters and to proceed from there.  While we can
determine this Jacobian \textit{a posteriori} by employing a cunning idea of Dumitriu and
Edelman, we do not have a direct derivation.

The idea of Dumitriu and Edelman can be summarized as follows:

\begin{lemma}\label{L:Jacobian}
Suppose $\phi:\mathcal{O}_1\to\mathcal{O}_2$ is an $N$-fold cover of $\mathcal{O}_2$ by $\mathcal{O}_1$,
both of which are open subsets of $\Reals^n$.  If the measure $f(x)\,d^nx$ is the symmetric pull-back
of the measure $g(y)\,d^ny$, with both $f$ and $g$ positive, then the Jacobian of $\phi$ is given by
$$
|\phi'(x)|=\frac{N f(x)}{g\circ\phi\,(x)}
$$
for any $x\in\mathcal{O}_1$.
\end{lemma}

\begin{proof}
For every $x\in\mathcal{O}_1$,
$$
f(x)\,d^nx = \phi^*\bigl(\tfrac{1}{N}g(y)\,d^ny\bigr) = \tfrac{1}{N}[g\circ\phi](x)|\phi'(x)| \,d^nx,
$$
which proves the lemma.
\end{proof}

Proposition~\ref{ItU} allows us to apply this lemma to the current situation.
The map of $\theta,\mu$ to the Verblunsky parameters is an $n!$-fold cover and so we obtain
\begin{align}\label{JU}
\left|\frac{\partial(\alpha,\phi)}{\partial(\theta,\mu)}\right|
&= 2^{1-n} \frac{|\Delta(e^{i\theta_1},\ldots,e^{i\theta_n})|^2 }
    {\prod_{k=0}^{n-2} (1-|\alpha_k|^2)^{n-k-2}}  \\
&= 2^{1-n} \frac{\prod_{k=0}^{n-2} (1-|\alpha_k|^2) }
    {\prod_{j=0}^{n} \mu_j } 
\end{align}
where one should regard the Verblunsky parameters as functions of the $\theta$s and $\mu$s.
The formulae above correspond to applying Proposition~\ref{ItU} with $\beta=2$ and $\beta=0$,
respectively.  Of course, one can also use other values of $\beta$, but the resulting formulae
are related to one another by Lemma~\ref{VandUn}.

Earlier we promised to determine the (common) normalization coefficient for the measures
\eqref{qqqqq} and \eqref{ggggg}.  We also promised to give a new derivation of the partition
function, \eqref{CGpart}, for the Coulomb gas. We will now settle these obligations.

\begin{lemma}\label{Somewhere}
The integral of \eqref{ggggg} is
\begin{equation}\label{SW1}
\int\!\!\cdots\!\!\int   \prod_{k=0}^{n-2}
(1-|\alpha_k|^2)^{\frac{\beta}{2}(n-k-1)-1} \,
    d^2 \alpha_0 \cdots d^2\alpha_{n-2}\,d\phi =
\frac{(2\pi)^n}{\beta^{n-1}(n-1)!}
\end{equation}
while
\begin{align}\label{SW2}
\int\!\! \cdots \!\! \int
\bigl|\Delta(e^{i\theta_1},\ldots,e^{i\theta_n})\bigr|^\beta \,
    \frac{d\theta_1}{2\pi} \cdots \frac{d\theta_n}{2\pi}
= \frac{\Gamma(\tfrac12\beta n + 1)}{\bigl[\Gamma(\tfrac12\beta +
1)\bigr]^n},
\end{align}
in agreement with \eqref{CGbeta} and \eqref{CGpart}.
\end{lemma}

\begin{proof}
Each of the integrals in \eqref{SW1} is rendered trivial by switching to polar coordinates:
$$
\int (1-|z|^2)^{(t/2)-1} \,d^2 z  = 2\pi t^{-1}.
$$
It is this integral that gives the normalization coefficient in the definition of the 
$\Theta$ distributions; cf. \eqref{E:ThetaDefn}.

The proof of \eqref{SW2} begins with the evaluation of the Dirichlet integral
$$
\int_\triangle \prod_{j=1}^n \mu_j^{\frac{\beta}{2}-1} \, d\mu_1 \cdots
d\mu_{n-1}
    = \frac{\Gamma(\tfrac{\beta}{2})^n}{\Gamma(\tfrac{n\beta}{2})},
$$
which is derived in the proof of Lemma~\ref{L:A4}; see
\eqref{DirichletI}.
As the two measures in Proposition~\ref{ItU} have the same total mass,
the integral of
\eqref{qqqqq} is given by \eqref{SW1}.  Combining this with the
Dirichlet integral above
leads us to
\begin{align*}
\text{LHS \eqref{SW2}}
=\frac{n!}{2^{1-n}}
\frac{1}{\beta^{n-1}(n-1)!}\frac{\Gamma(\tfrac{n\beta}{2})}{\Gamma(\tfrac{\beta}{2})^n}
=\frac{\Gamma(\tfrac{n\beta}{2}+1)}{\Gamma(\tfrac{\beta}{2}+1)^n},
\end{align*}
which is exactly \eqref{CGpart}.
\end{proof}

%%%%%%%%%%%%%%%%%%%%%%%%%%%%%%%%%%%%%%%%%%%%%%%%%%%%%%%%%%%%%%%%%%%%%%%%%%%%%%%%%%%%%%%%%%%%%%%%%%%%%%%%%
\section{The Proof of Theorem~2}
%%%%%%%%%%%%%%%%%%%%%%%%%%%%%%%%%%%%%%%%%%%%%%%%%%%%%%%%%%%%%%%%%%%%%%%%%%%%%%%%%%%%%%%%%%%%%%%%%%%%%%%%%

As explained in Section~\ref{S:2}, Theorem~\ref{T:2} is an immediate corollary of Proposition~\ref{OIt}
and the Geronimus relations.  As a result, the primary purpose of this section is to prove this proposition.

Throughout this section, $d\mu$ will denote a probability measure of the form given in \eqref{Ocoords}.
In particular, it is symmetric with respect to complex conjugation and the last Verblunsky parameter,
$\alpha_{2n-1}$, is equal to $-1$.   We will also use the notation $x_j=2\cos(\theta_j)$ repeatedly.

In addition to Lemma~\ref{VandUn} from the previous section, two further lemmas are required.
They are the following:

\begin{lemma}\label{VandVand}
If $x_j=2\cos(\theta_j)$, $1\leq j \leq n$, then
\begin{equation}
|\Delta(e^{\pm i\theta_1}, e^{\pm i\theta_2},\ldots, e^{\pm i\theta_n})|
    = |\Delta(x_1,x_2,\ldots,x_n)|^2 \, \prod_{l=1}^n
|2\sin(\theta_l)|
\end{equation}
where the left-hand side is shorthand for the Vandermonde of the $2n$ quantities
$e^{i\theta_1},e^{-i\theta_1}, \ldots, e^{i\theta_n}, e^{-i\theta_n}$.
\end{lemma}

\begin{proof}
By expanding the Vandermonde as in \eqref{VDefn},
\begin{align*}
&\phantom{{}={}} |\Delta(e^{\pm i\theta_1}, e^{\pm i\theta_2},\ldots, e^{\pm i\theta_n})| \\
&=\prod_{l=1}^n \bigl|e^{i\theta_l}-e^{-i\theta_l}\bigr|
    \prod_{j<k} |e^{i\theta_j}-e^{i\theta_k}|^2  |e^{i\theta_j}-e^{-i\theta_k}|^2 \\
&=\prod_{l=1}^n |2\sin(\theta_l)|
    \prod_{j<k} \bigl[2\cos(\theta_j) - 2\cos(\theta_k)\bigr]^2,
\end{align*}
as required.  In the last step we used that
$$
|(z-w)(z-\bar w)| = |(z-w)(1-\bar z \bar w)| = |(z+\bar z)-(w+\bar w)|
$$
for any pair of points $z,w$ on the unit circle.
\end{proof}

\begin{lemma}\label{PhiGamma}
Let $d\mu$ be a measure on the unit circle of the form given in \eqref{Ocoords}
and let $\Phi_k$ denote the corresponding monic orthogonal polynomials. Then
\begin{gather}
\prod_j (2-x_j) = \Phi_{2n}(1) = \prod_{k=0}^{2n-1} (1-\alpha_k) =
2\prod_{k=0}^{2n-2} (1-\alpha_k)    \label{Phi1}\\
\prod_j (2+x_j) = \Phi_{2n}(-1) = \prod_{k=0}^{2n-1} \bigl(1+(-1)^k
\alpha_k\bigr)              \label{Phi-1}
    = 2 \prod_{k=0}^{2n-2} \bigl(1+(-1)^k \alpha_k\bigr)
\end{gather}
where $x_j=2\cos(\theta_j)$ and $\alpha_k$ denote the Verblunsky parameters of $d\mu$.
\end{lemma}

\begin{proof}
As $\Phi_{2n}$ is orthogonal to each of $1,z,\ldots,z^{2n-1}$, which forms a basis for $L^2(d\mu)$,
its zeros must be $e^{\pm i\theta_1},\ldots,e^{\pm i\theta_n}$. Consequently, for any $x\in\Reals$,
$$
\Phi_{2n}(x) = \prod_j |x-e^{i\theta_j}|^2 = \prod_j
\bigl[x^2-2x\cos(\theta_j)+1\bigr],
$$
which gives the first equality in each of \eqref{Phi1} and \eqref{Phi-1}.

Because all Verblunsky parameters are real, the coefficients of the orthogonal polynomials are also real.
This implies $\Phi_k^*(z)=z^k\Phi_k(z^{-1})$ and so, for $z=\pm1$, the recurrence equation \eqref{PhiRec}
becomes $\Phi_{k+1}(z) = (z-\alpha_k z^k)\Phi_k(z)$.  Each of the second equalities stated in the lemma now
follows by the obvious induction.  The third equalities simply express $\alpha_{2n-1}=-1$.
\end{proof}

The following proposition is the $\SO(2n)$ analogue of Proposition~\ref{ItU} and so the main ingredient of
the proof of Theorem~\ref{T:2}.

\begin{prop}\label{OIt}
Consider the following measure on $[-2,2]^n\times\triangle:$
\begin{align*}
\frac{2^{-\kappa}}{n!} |\Delta(x_1,\ldots,x_n)|^{\beta}
\prod_{j=1}^n\mu_j^{\frac{\beta}{2}-1}
    \prod_{j=1}^n \bigl[(2-x_j)^a (2+x_j)^b\bigr] \; dx_1\cdots dx_n \, d\mu_1\cdots d\mu_{n-1}
\end{align*}
where $\kappa=(n-1)\frac{\beta}{2}+a+b+1$.  Under \eqref{Ocoords} and the change of variables
$x_j=2\cos\theta_j$, this gives a measure on $d\mu$ {\rm(}which is not normalized{\rm)}.
Transferring this measure to the Verbunsky parameters gives
\begin{align*}
    \prod_{k=0}^{2n-2} (1-\alpha_k^2)^{\frac{\beta(2n-k-1)}{4}-1}
    \prod_{k=0}^{2n-2} \bigl(1-\alpha_k\bigr)^{a+1-\frac\beta4}
             \bigl(1+(-1)^k\alpha_k\bigr)^{b+1-\frac\beta4} \; d\alpha_0 \cdots d\alpha_{2n-2}
\end{align*}
and $\alpha_{2n-1}\equiv 1$. After normalization, this measure corresponds to choosing the Verbunsky
parameters independently with distribution given by
\begin{equation}
\alpha_k \sim \begin{cases}
    B(\tfrac{2n-k-2}{4}\beta + a + 1,\tfrac{2n-k-2}{4}\beta + b + 1)& \text{$k$ even,} \\
    B(\tfrac{2n-k-3}{4}\beta + a+b+2,\tfrac{2n-k-1}{4}\beta)        & \text{$k$ odd.}
\end{cases}
\end{equation}
The definition of $B(s,t)$ is given in \eqref{E:beta}.
\end{prop}

\begin{proof}
From Lemma~\ref{VandVand} and Lemma~\ref{VandUn}, we may deduce that
\begin{equation}\label{E:Delta}
\begin{aligned}
|\Delta(x_1,x_2,\ldots,x_n)|^2 \prod_{l=1}^n |2\sin(\theta_l)|
    &= |\Delta(e^{\pm i\theta_1}, e^{\pm i\theta_2},\ldots, e^{\pm
i\theta_n})| \\
    &= 2^{n} \prod_{k=0}^{2n-2} (1-\alpha_k^2)^{(2n-k-1)/2} \prod_{j=1}^n
\mu_j^{-1}.
\end{aligned}
\end{equation}
By forming the square-root of the product of \eqref{Phi1} and
\eqref{Phi-1}, we can rewrite the product
of $|2\sin\theta_l|$ in terms of the Verblunsky parameters:
\begin{equation}\label{E:sin}
  \prod_{l=1}^n |2\sin(\theta_l)|
= 2 \prod_{k=0}^{2n-2} \bigl(1-\alpha_k\bigr)^{\frac12} \bigl(1+(-1)^k
\alpha_k\bigr)^{\frac12}.
\end{equation}
Substituting this into \eqref{E:Delta} above and doing a little
rearranging of terms leads us to
\begin{equation}\label{E:bbb}
2^{1-n} |\Delta(x_1,x_2,\ldots,x_n)|^2 \prod_{j=1}^n \mu_j
    = \frac{  \prod_{k=0}^{2n-2}
\bigl(1-\alpha_k^2\bigr)^{\frac{2n-k-1}{2}}  }%
      {  \prod_{k=0}^{2n-2} \bigl(1-\alpha_k\bigr)^{\frac12}
\bigl(1+(-1)^k \alpha_k\bigr)^{\frac12}  }.
\end{equation}
We will return to this equation in a moment.

For $\beta=2$ and $a=b=-\frac{1}{2}$, the proposition is an immediate corollary of
Propositions~\ref{SO2n} and~\ref{SO2nG} together with \eqref{E:sin}. The latter arises as the
Jacobian of the change of variables from $\theta_j$ to $x_j=2\cos(\theta_j)$.

Changing $a$ and $b$ amounts to multiplying this result by the appropriate powers of \eqref{Phi1}
and \eqref{Phi-1}, respectively.  To see that the two measures are equivalent for $\beta\neq2$,
it suffices to multiply by \eqref{E:bbb} raised to the $\frac{\beta}{2}-1$ power.
\end{proof}

Combining this proposition with Lemma~\ref{L:Jacobian} permits us to determine the
Jacobian of the map from the $(\theta,\mu)$-coordinates to the Verblunsky parameters.
From $\beta=2$, $a=b=-\frac{1}{2}$ we obtain
\begin{equation}\label{JO}
\left|\frac{\partial(\alpha)}{\partial(\theta,\mu)}\right| =
    \frac{2^{1-n}|\Delta(x_1,\ldots,x_n)|^2}%
    { \prod_{k=0}^{2n-2} (1-\alpha_k^2)^{\frac{2n-k-3}{2}} }.
\end{equation}
We do not have a direct derivation of this fact.

%%%%%%%%%%%%%%%%%%%%%%%%%%%%%%%%%%%%%%%%%%%%%%%%%%%%%%%%%%%%%%%%%%%%%%%%%%%%%%%%%%%%%%%%%%%%%%%%%%%%%%%%%
\section{The Selberg and Aomoto Integrals}\label{SandA}
%%%%%%%%%%%%%%%%%%%%%%%%%%%%%%%%%%%%%%%%%%%%%%%%%%%%%%%%%%%%%%%%%%%%%%%%%%%%%%%%%%%%%%%%%%%%%%%%%%%%%%%%%

In \cite{Selberg}, Selberg evaluated the following integral:
\begin{equation}\label{SI}
\int_0^1\!\!\cdots\!\!\int_0^1 |\Delta(u_1,\ldots,u_n)|^{2z}
\prod_{j=1}^n u_j^{x-1}(1-u_j)^{y-1}\,
    du_1\cdots du_n,
\end{equation}
which subsequently turned out to be important in random matrix theory.
We will present a new derivation of his result in a manner analogous to
the proof of Lemma~\ref{Somewhere} above.

We begin with the evaluation of the (common) integral of the measures given in Proposition~\ref{OIt}.
In the case of the second measure, this gives rise to a product of beta integrals (cf. \eqref{E:beta})
from which we obtain the answer
\begin{equation}\label{barf}
\begin{aligned}
{}& \prod_{\substack{ k=0 \\ \text{$k$ even} }}^{2n-2}
    \frac{  \Gamma(\frac{2n-k-2}{4}\beta+a+1)
\Gamma(\frac{2n-k-2}{4}\beta+b+1)  }%
        {  \Gamma(\frac{2n-k-2}{2}\beta+a+b+2)  }
2^{(2n-k-2)\frac{\beta}{2}+a+b+1} \\
{}\times{}& \prod_{\substack{ k=1 \\ \text{$k$ odd}  }}^{2n-3}
    \frac{  \Gamma(\frac{2n-k-3}{4}\beta+a+b+2)
\Gamma(\frac{2n-k-1}{4}\beta)  }%
        {  \Gamma(\frac{2n-k-2}{2}\beta+a+b+2)  }
2^{(2n-k-2)\frac{\beta}{2}+a+b+1}.
\end{aligned}
\end{equation}
To aid in the eventual comparison with Selberg \cite{Selberg}, we will switch to his parameters:
$$
x=a+1,\quad y=b+1, \quad\text{and}\quad z=\tfrac{1}{2}\beta.
$$
We also wish to make the following simplifications:  The products of the powers of $2$ can be
combined since they are the same for odd and even $k$, where we can easily sum the resulting
arithmetic progression in the exponents.  Similarly, we combine the products of the denominators
into a single product and make the substitution $p=2n-k-2$. In the even-$k$ numerators, we will
make the substitution $r=\frac12(2n-k-2)$ and in the odd-$k$ numerators,
the substitution $s=\frac12(2n-k-3)$.  Combining these we find that \eqref{barf} is equal to
\begin{align}\label{barf2}
2^\sigma \times
\frac{
    \prod_{r=0}^{n-1} \Gamma(rz+x)  \Gamma(rz+y)  \
    \prod_{s=0}^{n-2} \Gamma(sz+x+y)  \Gamma( (s+1)z )
}{
    \prod_{p=0}^{2n-2} \Gamma(pz+x+y)
}
\end{align}
where $\sigma=[(n-1)z+x+y-1](2n-1)$.  Note that there is a cancellation between the
$\Gamma(sz+x+y)$ terms in the numerator and the corresponding terms in the denominator.

This essentially completes the determination of the total mass of the measures in Proposition~\ref{OIt}.
After the transformation $x_j=4u_j-2$, the former of these is the tensor product of the
measure in \eqref{SI} with a measure in the $\mu_j$ coordinates.  This leads us to evaluate
\begin{align}\label{barf3}
 \frac{2^{-\kappa}}{n!} \int\!\!\cdots\!\!\int \prod_{j=1}^n\mu_j^{z-1} \,
d\mu_1\cdots d\mu_{n-1}
=\frac{2^{-\kappa}}{n!} \frac{\Gamma(z)^n}{\Gamma(nz)}
= \frac{2^{-\kappa} z^{1-n} \Gamma(z+1)^n}{(n-1)! \,\Gamma(nz+1)}
\end{align}
where $\kappa=(n-1)z+x+y-1$.  Moreover, taking the ratio \eqref{barf2}/\eqref{barf3} and
making the cancellation mentioned above, we obtain
\begin{equation}\label{SII}
\begin{aligned}
{}&\int\!\!\cdots\!\!\int |\Delta(x_1,\ldots,x_n)|^{2z}
    \prod_{j=1}^n \bigl[(2-x_j)^{x-1} (2+x_j)^{y-1}\bigr] \; dx_1\cdots
dx_n \\
{}={}&  2^{\sigma+\kappa} \times
\frac{
    \prod_{r=0}^{n-1} \Gamma(rz+x)  \Gamma(rz+y)  \  \prod_{s=0}^{n-2}
(s+1)z \Gamma( (s+1)z )
}{
    \Gamma(z+1)^n \, \prod_{p=n-1}^{2n-2} \Gamma(pz+x+y)
} \Gamma(nz+1).
\end{aligned}
\end{equation}
Notice that the term $z^{n-1}(n-1)!$ from \eqref{barf3} was split up inside the product over $s$.
To simplify, we use $\xi\Gamma(\xi)=\Gamma(\xi+1)$ inside the product over $s$ and notice that
the final factor, $\Gamma(nz+1)$, just corresponds to $s=n-1$.  In this way, all products run over
the same number of terms and can be combined.  Therefore, we reach the final conclusion
$$
\text{LHS \eqref{SII}} = 2^\tau \prod_{r=0}^{n-1} \frac{
\Gamma\bigl(rz+x\bigr) \Gamma\bigl(rz+y\bigr)
    \Gamma\bigl( (r+1)z \bigr)  }{ \Gamma\bigl(z+1\bigr)
\Gamma\bigl((n+r-1)z+x+y\bigr) }
$$
where $\tau=\sigma+\kappa=2n[(n-1)z+x+y-1]$.  This is in perfect agreement with Selberg's paper:
the values of \eqref{SI} and \eqref{SII} differ by a factor of $2^\tau$ as a result of the
change of variables $x_j=4u_j-2$; the measure is homogeneous of order $\tau/2$.

We now turn to our second topic: the Aomoto integral.

In \cite{Aomoto}, Aomoto determined the average value of $\prod(x-x_j)$ when the points
$x_j$ are distributed according to the Jacobi ensemble, \eqref{E:JE}.  Theorem~\ref{T:2}
shows that this is equivalent to evaluating the average of the characteristic polynomial
for a certain ensemble of Jacobi matrices.  In this way, Proposition~\ref{Aomo} below
reproduces Aomoto's result.

The answer is given in terms of the classical Jacobi polynomials:  In the notation of \cite{AbSt},
$$
\frac{4^n n!}{(a+b+n+1)_n} P_n^{(a,b)}(\tfrac{1}{2}x)
$$
are the monic polynomials that are orthogonal with respect to the measure
$$
\int f\,d\mu = \frac{\Gamma(a+b+2)}{4^{a+b+1}\Gamma(a+1)\Gamma(b+1)}
    \int_{-2}^2 f(x) \, (2-x)^a(2+x)^b \, dx.
$$
Here $(z)_n=\Gamma(z+n)/\Gamma(z)$ is the Pochhammer symbol. The recurrence coefficients for
the corresponding system of orthonormal polynomials are (cf. \eqref{OPRLrec})
\begin{equation}\label{JacRec}
\begin{aligned}
b_{n+1} &= \frac{2(b^2-a^2)}{(2n+a+b)(2n+a+b+2)} \\
a_{n+1}^2 &= \frac{16(n+1)(n+a+b+1)(n+a+1)(n+b+1)}{(2n+a+b+1)(2n+a+b+2)^2(2n+a+b+3)}.
\end{aligned}
\end{equation}

\begin{prop}\label{Aomo}
Let $J$ be the Jacobi matrix whose entries $b_1,a_1,\ldots,b_n$ are distributed as in Theorem~\ref{T:2}.
Then
\begin{equation}\label{Aom1}
    \mathbb{E}\bigl[\det(x-J)\bigr]
    = \frac{4^n n!}{(\tilde a+\tilde b+n+1)_n} P_n^{(\tilde a,\tilde b)}(\tfrac{1}{2}x)
\end{equation}
where $\tilde a=\frac{2(a+1)}{\beta}-1$ and $\tilde b=\frac{2(b+1)}{\beta}-1$.
\end{prop}

\begin{proof}
We will show that both sides of \eqref{Aom1} are related to the same monic
orthogonal polynomial on the unit circle.  We begin with the left-hand side.

Let $\alpha_0,\ldots,\alpha_{2n-2}$ be distributed as in \eqref{DistAR}, let $d\mu$
denote the corresponding measure on $S^1$, and let $d\nu$ denote the measure on $[-2,2]$
induced by \eqref{nuDefn}.  

The characteristic polynomial of $J$ is equal to $P_n(x)$, the $n$th monic orthogonal
polynomial associated to the measure $d\nu$; indeed the principal $k\times k$ minor
of $x-J$ is equal to $P_k$.  (By expanding along the last row, one can see that these
minor determinants obey the same recurrence relation as the orthogonal polynomials.)
Combining this with \eqref{Sz} gives
\begin{equation}\label{charJ}
\det((z+z^{-1})-J)=\frac{z^{-n}\Phi_{2n}(z)+z^n\Phi_{2n}(z^{-1})}{2}.
\end{equation}

As $\Phi_k$ is independent of $\alpha_k$ (it depends only on
$\alpha_0,\ldots,\alpha_{k-1}$), the recurrence relations \eqref{PhiRec} and \eqref{Phi*Rec} yield
\begin{align*}
\Exp \Phi_{k+1}(z)   &= z \, \Exp \Phi_k(z) - \Exp \bar{\alpha}_k \, \Exp \Phi_k^*(z)        \\
\Exp \Phi_{k+1}^*(z) &= \Exp \Phi_k^*(z) - z \, \Exp \alpha_k \, \Exp \Phi_k(z).
\end{align*}
Hence, $\Exp\Phi_{k}$ are the monic orthogonal polynomials associated to the averaged Verblunsky
parameters:
$$
\Exp(\alpha_k) = \begin{cases}
\frac{2b-2a}{(2n-k-2)\beta + 2a + 2b + 4}   & \text{$k$ even,} \\
\frac{\beta - 2a - 2b - 4}{(2n-k-2)\beta + 2a + 2b + 4}   & \text{$k$ odd.}
\end{cases}
$$
(Note that if $X\sim B(s,t)$ then $\Exp(X) = \frac{t-s}{t+s}$.)

Applying the Geronimus relations to these averaged Verblunsky coefficients does not produce 
the recursion coefficients associated to the Jacobi polynomials. However, by Proposition~\ref{P:B2},
$\Exp\Phi_{2n}$ is also the monic orthogonal polynomial of degree $2n$ associated to the `reversed' 
coefficients: $\tilde \alpha_k=\Exp\alpha_{2n-2-k}$, $0\leq k\leq 2n-2$, and $\tilde \alpha_{2n-1}=-1$.
Under the Geronimus relations, these coefficients give rise to
\begin{align*}
\tilde b_{k+1} &= \frac{2(\tilde b^2-\tilde a^2)}{(2k+\tilde a+\tilde b)(2k+\tilde a+\tilde b+2)} \\
\tilde a_{k+1}^2 &= \frac{16(k+1)(k+\tilde a+\tilde b+1)(k+\tilde a+1)
(k+\tilde b+1)}{(2k+\tilde a+\tilde b+1)(2k+\tilde a+\tilde b+2)^2(2k+\tilde a+\tilde b+3)}
\end{align*}
where $\tilde a=\frac{2(a+1)}{\beta}-1$ and $\tilde b=\frac{2(b+1)}{\beta}-1$. By comparison 
with \eqref{JacRec}, we find 
$$
\tfrac{1}{2}\Exp\{z^{-n}\Phi_{2n}(z)+z^n\Phi_{2n}(z^{-1})\}=
\tfrac{4^n n!}{(\tilde a+\tilde b+n+1)_n} P_n^{(\tilde a,\tilde b)}(\tfrac{1}{2}x).
$$ 
Equation \eqref{Aom1} now follows from \eqref{charJ}.
\end{proof}

\appendix

%%%%%%%%%%%%%%%%%%%%%%%%%%%%%%%%%%%%%%%%%%%%%%%%%%%%%%%%%%%%%%%%%%%%%%%%%%%%%%%%%%%%%%%%%%%%%%%%%%%%%%%%%
\section{The Surface Measure on $S^n$}
%%%%%%%%%%%%%%%%%%%%%%%%%%%%%%%%%%%%%%%%%%%%%%%%%%%%%%%%%%%%%%%%%%%%%%%%%%%%%%%%%%%%%%%%%%%%%%%%%%%%%%%%%

This Appendix presents some elementary results used in the text and is provided solely for the reader's convenience.

Let $d\tilde\sigma$ denote the usual surface measure on $S^n$ and let
$d\sigma$
denote the corresponding normalized probability measure on $S^n$.  We
will write
$\Disk^n$ for the $n$-disk: $\Disk^n=\{x\in\Reals^n : |x| < 1 \}$.

\begin{lemma}\label{L:A1}
If $f:\Reals^{n+1}\to\Cmplx$, then
$$
\int_{S^n} f(x) \, d\tilde\sigma(x) = \sum_{\pm} \int_{\Disk^n}
    f\bigl(x_1,\ldots,x_n,\pm\sqrt{1-|x|^2}\bigr) \frac{d^n
x}{\sqrt{1-|x|^2}}.
$$
\end{lemma}

\begin{proof}
The map $\phi:(x_1,\ldots,x_n)\mapsto (x_1,\ldots,x_n,\sqrt{1-|x|^2})$
is a diffeomorphism onto the
open upper hemisphere.  The corresponding Gram matrix is
$$
G_{i,j}
  = \Bigl\langle \frac{\partial \phi}{\partial x_i} \Big|
\frac{\partial \phi}{\partial x_j} \Bigr\rangle
  = \delta_{i,j} +
\frac{|x|^2}{1-|x|^2}\frac{x_i}{|x|}\frac{x_j}{|x|},
$$
which is a rank one perturbation of the identity matrix.  Therefore,
$$
\det[G]=1+\frac{|x|^2}{1-|x|^2} = \frac{1}{1-|x|^2},
$$
from which the lemma follows immediately.
\end{proof}

We wish to determine the distributions of certain probability measures
induced from the normalized surface
measure on $S^n$. These follow easily from this lemma.

\begin{coro}\label{C:A2}
For the projection of the $n$-sphere into its first coordinate, we
have
\begin{align}
\int_{S^n} f(x_1) \, d\sigma(x) =
\frac{2^{1-n}\Gamma(n)}{\Gamma(\frac{n}{2})^2}
    \int_{-1}^1 f(t) \, (1-t^2)^{(n-2)/2} \, dt,
\end{align}
that is, $x_1$ is $B(\tfrac{n}{2},\tfrac{n}{2})$-distributed.
Projection onto the first two coordinates gives
\begin{align}
\int_{S^n} f(x_1 + i x_2) \, d\sigma(x) =  \frac{n-1}{2\pi}
    \int_{\Disk} f(z) \, (1-|z|^2)^{(n-3)/2} \, d^2z,
\end{align}
which implies $x_1+ix_2$ is $\Theta_{n}$-distributed.
\end{coro}

\begin{proof}  Both formulae follow from the more general statement
that for $1\leq k < n$,
$$
\int_{S^n} f(x_1,\ldots,x_k) \, d\sigma(x) \propto
    \int_{\Disk^k} f(x_1,\ldots,x_k) (1-|x|^2)^{(n-k-1)/2} \, dx_1
\cdots dx_k
$$
where the proportionality constant depends on $k$, but not the function $f$.  The value
of the normalization constant can then be determined by substituting $f \equiv 1$.

An inductive proof of the more general statement follows easily from
$$
  \int_{-s}^s \bigl(s^2-x_{k+1}^2\bigr)^{(n-k-2)/2} \, dx_{k+1} \propto s^{n-k-1}
$$
where one takes $s = (1-x_1^2-\cdots-x_k^2)^{1/2}$.
\end{proof}

\begin{coro}\label{C:A3}
If the vector $(x_1,y_1,\ldots,x_n,y_n)$ is chosen at random from the $(2n-1)$-sphere
according to normalized surface measure, then
$$
(\mu_1,\ldots,\mu_n) = \bigl(x_1^2+y_1^2,\ldots,x_n^2+y_n^2\bigr)
$$
is uniformly distributed on the $(n-1)$-simplex $\sum \mu_j =1$. That is,
$$
\Exp f(\mu_1,\ldots,\mu_n) = (n-1)!\int_{\triangle}
f(s_1,\ldots,s_{n-1},1-s_1-\cdots-s_{n-1}) \, ds_1\cdots ds_{n-1}
$$
where $\triangle=\{(s_1,\ldots,s_{n-1}) \, : \, 0\leq s_1+\cdots+s_{n-1}\leq 1\}$.
\end{coro}

\begin{proof}
By Lemma \ref{L:A1}, we need only compute
$$
\sum_{\pm} \int_{\Disk^{2n-1}} f(x_1,y_1,\ldots,x_n,\pm y_n) \, \frac{dx_1dy_1\cdots dx_n}{y_n}
$$
where $y_n=\sqrt{1-(x_1^2+y_1^2+\cdots +x_n^2)}$.

Let us change variables by
\begin{align*}
&\left.
    \begin{aligned} x_j&=\cos(\phi_j)\sqrt{s_j} \\
    y_j&=\sin(\phi_j)\sqrt{s_j} \end{aligned}
\quad \right\} 1\leq j\leq n-1 \\
&\left.x_n=\cos(\phi_n)\biggl(1-\sum_{k=1}^{n-1}s_k\biggr)^{1/2},\right.
\end{align*}
for which the Jacobian is $y_n^{-1} dx_1dy_1\cdots dx_n = ds_1 ds_2\cdots ds_{n-1}d\phi_1\ldots d\phi_n$.
Hence, up to proportionality constants that depend only on $n$,
$$
\Exp f(\mu_1,\ldots,\mu_n) \propto \int_{\triangle} f(s_1,\ldots,s_{n-1},1-s_1-\cdots-s_{n-1})
\, ds_1\cdots ds_{n-1}.
$$

The determination of the normalization constant is an easy calculation.
A more general result will appear in the proof of Lemma~\ref{L:A4}.
\end{proof}

\begin{lemma}\label{L:A4}
For $(\mu_1,\ldots,\mu_n)$ uniformly distributed on the $(n-1)$-simplex $\sum \mu_j =1$ and
$p_1,\ldots,p_n$ real numbers greater than $-1$, we have
\begin{align}
\Exp \bigl\{ \prod \mu_j^{p_j} \bigr\} = \frac{(n-1)!\,
\Gamma(p_1+1) \cdots \Gamma(p_n+1)}{\Gamma(p_1+\cdots+p_n+n)}.
\end{align}
\end{lemma}

\begin{proof}
This is sometimes known as Dirichlet's integral; our method of proof is a standard one.
We begin with the product integral
$$
\prod_{j=1}^{n} \int_0^{\infty} s_j^{p_j} e^{-s_j}\,ds_j =
\prod_{j=1}^{n}\Gamma(p_j+1)
$$
and change variables to $r=\sum s_j$ and $\mu_j=s_j/r$.  The Jacobian is given by
$r^{n-1}\,dr \, d\mu_1\cdots d\mu_{n-1}=ds_1\cdots ds_n$ which leads us to
\begin{align*}
\prod_{j=1}^{n}\Gamma(p_j+1) &=
     \int_{\triangle} \! \int_0^\infty \prod (r\mu_j)^{p_j} \, e^{-r}
r^{n-1} \, dr \, d\mu_1\cdots d\mu_{n-1} \\
&= \Gamma(p_1+\cdots+p_n+n) \,
    \int_{\triangle} \prod \mu_j^{p_j} \, d\mu_1\cdots d\mu_{n-1}.
\end{align*}

By taking $p_1=\cdots=p_n=0$ we recover the normalization constant of the $\mu_j$-integral given in Corollary~\ref{C:A3}.  Moreover,
\begin{equation}\label{DirichletI}
\begin{aligned}
  \Exp \bigl\{ \prod \mu_j^{p_j} \bigr\}
&= (n-1)! \int_{\triangle} \prod \mu_j^{p_j} \, d\mu_1\cdots d\mu_{n-1}
\\
&= \frac{(n-1)!\,\Gamma(p_1+1) \cdots
\Gamma(p_n+1)}{\Gamma(p_1+\cdots+p_n+n)}
\end{aligned}
\end{equation}
as promised.
\end{proof}

%%%%%%%%%%%%%%%%%%%%%%%%%%%%%%%%%%%%%%%%%%%%%%%%%%%%%%%%%%%%%%%%%%%%%%%%%%%%%%%%%%%%%%%%%
\section{The CMV Matrix and Geronimus Relations}\label{S:CMV}
%%%%%%%%%%%%%%%%%%%%%%%%%%%%%%%%%%%%%%%%%%%%%%%%%%%%%%%%%%%%%%%%%%%%%%%%%%%%%%%%%%%%%%%%%

The purpose of this appendix is to describe some results in the theory
of orthogonal polynomials on the unit circle that were used in the main text;
we also wish to show how easily they can be derived from the
perspective of Cantero, Moral, and Vel\'azquez, \cite{CMV}.
We begin with an outline of that work.

Applying the Gram--Schmidt procedure to
$1,z,z^{-1},z^2,z^{-2},\ldots$ in $L^2(d\mu)$ produces the orthonormal basis
\begin{equation}
\chi_k(z) = \begin{cases}
    z^{-k/2}    \phi_k^*(z) & \text{$k$ even} \\
    z^{(1-k)/2} \phi_k(z) & \text{$k$ odd}\end{cases}
\end{equation}
where $k\geq0$.  As in the Introduction, $\phi_k$ denotes the $k$th orthonormal polynomial and
$\phi_k^*$, its reversal (cf. \eqref{rev}).
If we apply the procedure to $1,z^{-1},z,z^{-2},z^2,\ldots$ instead, then we obtain
a second orthonormal basis:
\begin{equation}
X_k(z)=\overline{\chi_k(1/\bar z)}= \begin{cases}
    z^{-k/2}    \phi_k(z) & \text{$k$ even} \\
    z^{(-1-k)/2} \phi_k^*(z) & \text{$k$ odd.}\end{cases}
\end{equation}

It is natural to compute the matrix representation of $f(z)\mapsto zf(z)$ in $L^2(d\mu)$ with respect
to these bases.  This was done in a rather cunning way.  The matrices with entries
$$
L_{i+1,j+1} = \langle \chi_i(z)| z X_j(z) \rangle \quad\text{and}\quad
M_{i+1,j+1} = \langle   X_i(z) | \chi_j(z)\rangle
$$
are block-diagonal; indeed,
\begin{equation}\label{LMdefn}
L=\diag\bigl(\Xi_0   ,\Xi_2,\Xi_4,\ldots\bigr) \quad\text{and}\quad
M=\diag\bigl(\Xi_{-1},\Xi_1,\Xi_3,\ldots\bigr)
\end{equation}
where $\Xi_{-1}=[1]$ and
$$
\Xi_k = \begin{bmatrix} \bar\alpha_k & \rho_k \\ \rho_k & -\alpha_k
\end{bmatrix}.
$$
The representation of $f(z)\mapsto zf(z)$ in the $\chi_j$ basis is just
$LM$, which is the CMV matrix of the title of this appendix, and in the $X_j$ basis, it is $ML$.

If $d\mu$ is supported at finitely many points, say $m$, then $\alpha_{m-1}$ will be unimodular
and hence $\Xi_{m-1}$ is diagonal.  We replace $\Xi_{m-1}$ by the $1\times 1$ matrix that is its
top left entry, $\bar\alpha_{m-1}$, and discard all $\Xi_k$ with $k\geq m$.
In this way, we find $L$ and $M$ are naturally $m\times m$ block-diagonal matrices.

In a couple of places we noted that the spectral measure for any $U\in\SO(2n)$ has last Verblunsky
parameter equal to $-1$.  The CMV matrix allows us to give a particularly short proof of this fact:

\begin{lemma}
Let $d\mu$ be a probability measure of the form
$$
\int f \,d\mu = \sum_j \mu_j f(z_j)
$$
where $z_1,\ldots,z_m$ are distinct points on the unit circle.  If $\alpha_{m-1}$
denotes the final Verblunsky parameter associated to this measure, then
\begin{equation}\label{detCMV}
\prod_{j=1}^m z_j = (-1)^{m-1} \bar\alpha_{m-1}.
\end{equation}
\end{lemma}

\begin{proof}
As one might guess, \eqref{detCMV} just represents two ways of calculating the determinant
of the CMV matrix: the left-hand side is the product of the eigenvalues; the right-hand side
is the product of the determinants of $L$ and $M$.  Note that $\det(\Xi_k)=-1$ for $0 \leq k < m$
while $\det(\Xi_{m-1})=\bar\alpha_{m-1}$.
\end{proof}

Next we prove a result used in the derivation of the Aomoto integral:

\begin{prop}\label{P:B2}
Given a finite system of Veblunsky parameters $\alpha_k\in\Disk$,
$0\leq k\leq m-2$, and $\alpha_{m-1}=e^{i\phi}$, define a
second system by $\tilde\alpha_k=-e^{i\phi}\bar\alpha_{m-2-k}$, $0\leq k \leq m-2$, and
$\tilde\alpha_{m-1}=e^{i\phi}$. Then the monic orthogonal
polynomials of degree $m$ associated to these two systems are the same.
\end{prop}

\begin{proof}
If $L$ and $M$ are the matrices associated to the $\alpha$s and $\Phi_m$ the monic orthogonal
polynomial of degree $m$, then $\det(z-LM)=\Phi_m(z)$.  This is because both are monic polynomials
vanishing at the eigenvalues of $LM$.

Similarly, we have  $\det(z-\tilde L \tilde M)=\tilde \Phi_m(z)$ for the $\tilde\alpha$s.  Hence
the proposition will follow once we show that $LM$ and $\tilde L \tilde M$ are conjugate.  We will
give full details when $m$ is even and a few remarks on the changes necessary when $m$ is odd.

Conjugating $LM$ by
$$
\begin{bmatrix}
 0 & \ldots  &   0   &   1  &   \\
 0 & \ldots  &   1   &   0  &   \\
 \vdots&\ddots&\vdots&\vdots&   \\
 1 &\ldots   &   0   &   0  &
\end{bmatrix}
$$
is equivalent to reversing the order of the rows and columns in each factor.
As $m$ is even, both $L$ and $M$ retain their structure. The transformation on
$L$ amounts to replacing $\alpha_k$ by $-\bar\alpha_{m-2-k}$; for $M$, we have
the additional change that the $(1,1)$ entry becomes $e^{-i\phi}$ and the $(m,m)$ entry, $1$.

To convert $M$ to its proper structure (with $1$ in the first position), we perform the
transformations $L\mapsto U^\dagger L V$ and $M\mapsto V^\dagger M U$ where
$$
U=\diag(e^{i\phi},1,e^{i\phi},1,\ldots), \quad V=\diag(1,e^{i\phi},1,e^{i\phi},\ldots).
$$
In this way, we obtain the matrices $\tilde L$ and $\tilde M$, which shows that
$\tilde L\tilde M$ is conjugate to the original $LM$.

When $m$ is odd, reversing the order of the rows and columns converts $L$ to a matrix
whose structure resembles that of $M$, while $M$ is converted to an $L$-like matrix.
Proceeding as above shows that $LM$ is conjugate to $\tilde M\tilde L$, and hence to $\tilde L\tilde M$.
\end{proof}

Let us now consider the case where the measure $d\mu$ is symmetric with respect to complex conjugation,
or what is equivalent, where all Verblunsky parameters are real.  It is a famous observation of Szeg\H{o}
(see \cite[\S11.5]{Szego}) that the polynomials orthogonal with respect to this measure are intimately
related to the polynomials orthogonal with respect to the measure $d\nu$ on $[-2,2]$ defined by
\begin{equation}\label{nuDefn}
  \int_{S^1} f(z+z^{-1}) \, d\mu(z)  =  \int_{-2}^2 f(x) \,d\nu(x).
\end{equation}

As described in Section~\ref{S:2}, the recurrence coefficients for these systems of orthogonal
polynomials are related by Geronimus relations:
\begin{equation}\label{AB:Geron}
\begin{aligned}
b_{k+1}   &= (1-\alpha_{2k-1})\alpha_{2k} -
(1+\alpha_{2k-1})\alpha_{2k-2}           \\
a_{k+1}   &= \big\{ (1-\alpha_{2k-1})(1-\alpha_{2k}^2)(1+\alpha_{2k+1})
\big\}^{1/2}.
\end{aligned}
\end{equation}
We will now present a short proof of these formulae.  As an off-shoot of our method, we also
recover relations to the recurrence coefficients for $(4-x^2)\,d\nu(x)$ and $(2\pm x)d\nu$.
The former appears in the proposition below, the latter in the remark that follows it.

\begin{prop} Let $\alpha_k$ be the system of real Verblunsky parameters associated to a
symmetric measure $d\mu$ and let $L$ and $M$ denote the matrices of \eqref{LMdefn}.
Then $LM+ML$ is unitarily equivalent to the direct sum of two Jacobi matrices:
$$
J = \begin{bmatrix}
 b_1 & a_1  &   0   &      \\
 a_1 & b_2  &\ddots &      \\
  0  &\ddots&\ddots &      \\
\end{bmatrix}
\qquad
\tilde J = \begin{bmatrix}
 \tilde b_1 & \tilde a_1  &   0   &      \\
 \tilde a_1 & \tilde b_2  &\ddots &      \\
     0      &\ddots       &\ddots &      \\
\end{bmatrix}
$$
where $a_k$ and $b_k$ are as in \eqref{AB:Geron} and
\begin{align*}
\tilde b_{k+1}   &= (1-\alpha_{2k+1})\alpha_{2k} - (1+\alpha_{2k+1})\alpha_{2k+2}           \\
\tilde a_{k+1}   &= \big\{ (1+\alpha_{2k+1})(1-\alpha_{2k+2}^2)(1-\alpha_{2k+3}) \big\}^{1/2}.
\end{align*}
Moreover, the spectral measure for $(J,e_1)$ is precisely the $d\nu$ of \eqref{nuDefn}.
The spectral measure for $(\tilde J,e_1)$ is $\frac1{2(1-\alpha_0^2)(1-\alpha_1)}(4-x^2)\,d\nu(x)$.
\end{prop}

\begin{proof}
Let $S$ denote the following unitary block matrix
$$
S=\diag([1],S_1,S_3,\ldots) \quad\text{where}\quad
    S_{k}= \frac{1}{\sqrt{2}}
    \begin{bmatrix} -\sqrt{1-\alpha_k} & \sqrt{1+\alpha_k} \\
        \sqrt{1+\alpha_k} & \sqrt{1-\alpha_k}\end{bmatrix},
$$
which is easily seen to diagonalize $M$.  Indeed, $S^\dagger M S =\diag(+1,-1,+1,-1,\ldots)$.
We will denote this matrix by $R$.

The matrix $LM+ML$ is unitarily equivalent to $A=S^\dagger(LM+ML)S=S^\dagger LSR+ RS^\dagger LS$,
which we will show is the direct sum of two Jacobi matrices.  We begin by showing that even-odd
and odd-even entries of $A$ vanish, from which it follows that $A$ is the direct sum of its even-even
and odd-odd submatrices.

Left multiplication by $R$ changes the sign of the entries in each even-numbered row, while right
multiplication by $R$ reverses the sign of each even-numbered column. In this way, $RB+BR$ has
the stated direct sum structure for any matrix $B$ and hence, in particular, for $B=S^\dagger LS$.

It remains only to calculate the non-zero entries of $A$.  As $S$ and $L$ are both tri-diagonal,
$A$ must be hepta-diagonal and so the direct sum of tri-diagonal matrices.  Moreover, $A$
is symmetric (because $L$ is) so there are only four categories of entries to calculate:
the odd/even diagonals and the odd/even off-diagonals.  We begin with the diagonals:
\begin{align*}
A_{2k+1,2k+1} &= \begin{bmatrix} \sqrt{1+\alpha_{2k-1}} & \sqrt{1-\alpha_{2k-1}} \end{bmatrix}
    \begin{bmatrix} -\alpha_{2k-2} & 0 \\ 0 & \alpha_{2k} \end{bmatrix}
    \begin{bmatrix} \sqrt{1+\alpha_{2k-1}} \\ \sqrt{1-\alpha_{2k-1}} \end{bmatrix} \\
&= (1-\alpha_{2k-1})\alpha_{2k} - (1+\alpha_{2k-1})\alpha_{2k-2} \\
A_{2k,2k} &= -\begin{bmatrix} -\sqrt{1-\alpha_{2k-1}} & \sqrt{1+\alpha_{2k-1}} \end{bmatrix}
    \begin{bmatrix} -\alpha_{2k-2} & 0 \\ 0 & \alpha_{2k} \end{bmatrix}
    \begin{bmatrix} -\sqrt{1-\alpha_{2k-1}} \\ \sqrt{1+\alpha_{2k-1}} \end{bmatrix} \\
&= (1-\alpha_{2k-1})\alpha_{2k-2}-(1+\alpha_{2k-1})\alpha_{2k}.
\end{align*}
Note that the factor of 2 resulting from $A$ being the sum of two terms is
cancelled by the factors of $2^{-1/2}$ coming from $S$ and $S^\dagger$.
The calculation of the off-diagonal terms proceeds in a similar fashion:
\begin{align*}
A_{2k+1,2k+3} &= \begin{bmatrix} \sqrt{1+\alpha_{2k-1}} &
\sqrt{1-\alpha_{2k-1}} \end{bmatrix}
    \begin{bmatrix} 0 & 0 \\ \rho_{2k} & 0 \end{bmatrix}
    \begin{bmatrix} \sqrt{1+\alpha_{2k+1}} \\ \sqrt{1-\alpha_{2k+1}} \end{bmatrix} \\
&= \sqrt{(1-\alpha_{2k-1})(1-\alpha_{2k}^2)(1+\alpha_{2k+1})} \\
A_{2k,2k+2} &= -\begin{bmatrix} -\sqrt{1-\alpha_{2k-1}} & \sqrt{1+\alpha_{2k-1}} \end{bmatrix}
    \begin{bmatrix} 0 & 0 \\ \rho_{2k} & 0 \end{bmatrix}
    \begin{bmatrix} -\sqrt{1-\alpha_{2k+1}} \\ \sqrt{1+\alpha_{2k+1}} \end{bmatrix} \\
&= \sqrt{(1+\alpha_{2k-1})(1-\alpha_{2k}^2)(1-\alpha_{2k+1})}.
\end{align*}

That $d\nu$ is the spectral measure for $(J,e_1)$ is an immediate consequence
of the spectral theorem, $LM+ML=LM+(LM)^{-1}$, and the fact that $S$ leaves the vector
$[1,0,\ldots,0]$ invariant.

Tracing back through the definitions, we find that the spectral measure for $(\tilde J,e_1)$
is equal to that for the operator $f(z)\mapsto (z+z^{-1})f(z)$ in $L^2(d\mu)$ and the vector
$$
f(z)=\left(\tfrac{1+\alpha_1}{2}\right)^{\frac12} \chi_2(z) -
    \left(\tfrac{1-\alpha_1}{2}\right)^{\frac12} \chi_1(z)
= \left(\tfrac{1+\alpha_1}{2}\right)^{\frac12} z^{-1} \phi_2^*(z) -
    \left(\tfrac{1-\alpha_1}{2}\right)^{\frac12} \phi_1(z).
$$
From the relations \eqref{PhiRec}, \eqref{Phi*Rec}, and \eqref{PhiNorm}, we find
$\rho_1\phi_2^*(z)=\phi_1^*(z) - \alpha_1 z\phi_1(z)$, $\rho_0\phi_1(z)=z-\alpha_0$,
and $\rho_0\phi_1^*(z)=1-\alpha_0z$.  These simplify the formula considerably:
$$
f(z) = \frac{z^{-1} - z }{\rho_0\sqrt{2(1-\alpha_1)}}.
$$
The expression for the spectral measure for $(\tilde J,e_1)$ now
follows from the simple
calculation $|z^{-1}-z|^2=4-(z+z^{-1})^2$.
\end{proof}

\noindent
\textit{Remark:} In the above proof, we conjugated $LM+ML$ by the unitary matrix which diagonalizes
$M$.  One may use the matrix $\diag(S_0,S_2,\ldots)$, which diagonalizes $L$, instead.  This also
conjugates $LM+ML$ to the direct sum of two Jacobi matrices.  In this way, we learn that
the recurrence coefficients for $\tfrac{1}{2(1\pm\alpha_0)}(2\pm x)d\nu(x)$ are given by 
\begin{align*}
b_{k+1} &= \pm(1\mp\alpha_{2k})\alpha_{2k+1} \mp (1\pm\alpha_{2k})\alpha_{2k-1} \\
a_{k+1} &= \big\{(1\mp\alpha_{2k})(1-\alpha_{2k+1}^2)(1\pm\alpha_{2k+2})\big\}^{1/2}.
\end{align*}

%%%%%%%%%%%%%%%%%%%%%%%%%%%%%%%%%%%%%%%%%%%%%%%%%%%%%%%%%%%%%%%%%%%%%%%%%%%%%%%%%%%%%%%%%

\end{document}